\documentclass[a4paper,reqno,12pt, twoside]{amsart}
\usepackage{amsfonts}
\usepackage{caption}
\usepackage{color}
\usepackage{hyperref}
\definecolor{ddorange}{rgb}{1,0.5,0}
\definecolor{ddcyan}{rgb}{0,0.2,1.0}

\usepackage[left=2.5cm,right=2.5cm,top=2.2cm,bottom=2.3cm,headsep=0.7cm]{geometry}

\usepackage[T1]{fontenc} 
\usepackage[utf8]{inputenc}
\usepackage[french, italian, english]{babel}
\usepackage{amsmath,amssymb,amsthm,mathrsfs,esint}
\usepackage{mathtools}

\newcommand{\R}{{\mathbb R}}
\newcommand{\N}{{\mathbb N}}
\newcommand{\Z}{{\mathbb Z}}
\newcommand{\E}{{\mathrm E}}
\newcommand{\m}{{\textbf{m}}}

\newcommand{\km}{{k^{-1}}}
\newcommand{\tq}{{\tilde{q}}}
\newcommand{\qz}{{q_z^k}}
\newcommand{\tqz}{{\tilde{q}_z^k}}
\newcommand{\Qz}{{Q_z^k}}
\newcommand{\Rn}{{\R}^n}
\newcommand{\dx}{\, \mathrm{d} x}
\renewcommand{\dh}{\, \mathrm{d} \mathcal{H}^{n-1}}
\newcommand{\hn}{\mathcal{H}^{n-1}}
\newcommand{\Sn}{{\mathbb{S}^{n-1}}}
\newcommand{\ol}{\overline}

\newcommand{\tr}{\mathrm{tr}}

\newcommand{\sm}{\setminus}
\newcommand{\Mnn}{{\mathbb{M}^{n\times n}_{sym}}}

\newcommand{\weak}{\rightharpoonup}
\newcommand{\wstar}{\stackrel{*}\rightharpoonup}

\newcommand{\mres}{\mathbin{\vrule height 1.6ex depth 0pt width
0.13ex\vrule height 0.13ex depth 0pt width 1.3ex}}

\newcommand{\interior}[1]{%
  {\kern0pt#1}^{\mathrm{o}}%
}

\theoremstyle{plain}
\begingroup
\theoremstyle{plain}
\newtheorem{theorem}{Theorem}[section]

\newtheorem{proposition}[theorem]{Proposition}
\newtheorem{lemma}[theorem]{Lemma}
\theoremstyle{definition}

\theoremstyle{remark}
\newtheorem{remark}[theorem]{Remark}

\endgroup

\numberwithin{equation}{section}
\title[On the approximation of $SBD$ functions and applications]{
On the approximation of $SBD$ functions \\ and some applications}
\author{Vito Crismale}
\address{CMAP, \'Ecole Polytechnique, UMR CNRS 7641, 91128 Palaiseau Cedex, France}
\email[Vito Crismale]{vito.crismale@polytechnique.edu}

\begin{document}
\begin{abstract}
Three density theorems for three suitable subspaces of $SBD$ functions, in the strong $BD$ topology, are proven. 
The spaces are $SBD$, $SBD^p_\infty$, where the absolutely continuous part of the symmetric gradient is in $L^p$, with $p>1$,
and $SBD^p$, whose functions are in $SBD^p_\infty$ and the jump set has finite $\hn$-measure. This generalises on the one hand the density result \cite{Cha04} by Chambolle and, on the other hand, extends in some sense the three approximation theorems in \cite{DPFusPra17} by De Philippis, Fusco, Pratelli for $SBV$, $SBV^p_\infty$, $SBV^p$ spaces, obtaining also more regularity for the absolutely continuous part of the approximating functions. As application, the sharp version of two $\Gamma$-convergence results for energies defined on $SBD^2$ is derived.
\end{abstract}
\maketitle

{\small

\tableofcontents

\section{Introduction}

The study of free discontinuity functionals has required the introduction of suitable ambient spaces, such as the Functions of Bounded Variations $BV$ and of Bounded Deformation $BD$, with corresponding subspaces and generalisations.

A $L^1$ function $u$ is in $BV$ [respectively in $BD$] if its distributional gradient $\mathrm{D}u$ [resp.\ its distributional symmetric gradient $\E u=(\mathrm{D}u+\mathrm{D}^Tu)/2$] is a bounded Radon measure. In particular, a $BD$ function is defined from a set $\Omega\subset\Rn$ into $\Rn$. 
The measure $\mathrm{D}u$ [$\E u$] is decomposed into three parts: one absolutely continuous with respect to $\mathcal{L}^n$, with density $\nabla u$ [$e(u)$], one supported on the rectifiable $(n{-}1)$-dimensional \emph{jump set} $J_u$, where $u$ has two different approximate limits $u^+$, $u^-$ on the two sides of $J_u$ with respect to an approximate normal $\nu_u\in \Sn$, and a \emph{Cantor part}, vanishing on Borel sets of finite $\hn$ measure. $SBV$ [$SBD$] is the space of $BV$ [$BD$] functions with null Cantor part.  
Here we consider also, for $p>1$, the subspaces
\begin{equation*}
SBD^p(\Omega):=\{u\in SBD(\Omega)\colon e(u)\in L^p(\Omega;\Mnn)\,,\ \hn(J_u)<\infty\}
\end{equation*}
and
\begin{equation*}
SBD^p_\infty(\Omega):=\{u\in SBD(\Omega)\colon e(u)\in L^p(\Omega;\Mnn)\}\,
\end{equation*}
with analogous definitions for $SBV^p(\Omega)$ and $SBV^p_\infty$ (see Section~\ref{Sec1} for more details). 

The spaces $SBD^p$ are very important in Fracture Mechanics: 
if $u$ represents the \emph{displacement} of a body from its equilibrium configuration, $J_u$ is nothing but the \emph{crack} set and $e(u)$ is the \emph{linearised elastic strain}, which is in $L^2$ (so $p=2$) if the material is \emph{linearly elastic} in the bulk region. For many years after the introduction of $SBD$ in \cite{AmbCosDM97}, $SBD^2$ has been employed to study \emph{brittle fracture}, namely the Griffith energy
\begin{equation}\tag{G}\label{eq:Griffith}
\int\limits_\Omega \mathbb{C}e(u)\colon e(u)\dx + \hn(J_u)\,,
\end{equation}
$\mathbb{C}$ being the (fourth-order positive definite) Cauchy stress tensor, with possibly lower order terms due to forces, and boundary conditions. Unfortunately, the corresponding compactness and lower semicontinuity theorem 
\cite{BelCosDM98} requires equi-integrability of displacements, which is not guaranteed for a sequence with bounded \eqref{eq:Griffith} energy. Indeed, the right ambient space for \eqref{eq:Griffith} is $GSBD^2$, introduced by Dal Maso in \cite{DM13}, with the corresponding compactness and lower semicontinuity theorem proven very recently in \cite{CC18} (see also \cite{FriSol16} in dimension 2).

The first density result for $SBD^2$, due to Chambolle (\cite{Cha04, Cha05Add}) consists then in the approximation, \emph{with respect to the energy} \eqref{eq:Griffith}, of $u\in SBD^2 \cap L^2$ by functions smooth outside their jump set, in turn closed and included in a finite union of $C^1$ hypersurfaces (this has been extended to $GSBD^p$ in \cite{FriPWKorn, CFI17Density, CC17}). 

If we are given an energy controlling the \emph{amplitude} of the jump $[u]:=u^+-u^-$ in $L^1(J_u;\Rn)$, in contrast to Griffith energy that controls only the measure of $J_u$, then $SBD^p$ (for a $p$-growing bulk energy) is the proper ambient space. This is the case, for $p=2$, of the energy
\begin{equation}\tag{C}\label{eq:cohesive}
\int\limits_\Omega \mathbb{C}e(u)\colon e(u)\dx + \hn(J_u) + \int \limits_{J_u} \big| [u] \odot \nu_u \big| \dh\,,
\end{equation}
($\odot$ being the symmetric tensor product) considered by Focardi and Iurlano in \cite{FocIur14}, and recently in \cite{CarVG18}. A fracture energy depending on $[u]$, as \eqref{eq:cohesive}, is often called \emph{cohesive}, in contrast to the brittle energy \eqref{eq:Griffith}.

In order to deal with energies such as \eqref{eq:cohesive}, the following approximation theorem for $SBD^p$, that involves also the jump part of $\E u$, is proven. 
This is 
 the main result of the paper.
\begin{theorem}\label{teo:density}
Let $\Omega$ be an open bounded Lipschitz subset of $\Rn$, and $u\in SBD^p(\Omega)$, with $p>1$.
Then there exist $u_k\in SBV^p(\Omega;\Rn)\cap L^\infty(\Omega; \Rn)$ such that each
$J_{u_k}$ is closed 
and included in a finite union of closed connected pieces of $C^1$ hypersurfaces, $u_k\in C^\infty(\ol \Omega\sm J_{u_k};\Rn) \cap W^{m,\infty}(\Omega\setminus J_{u_k}; \Rn)$ for every $m\in \N$, and:
\begin{subequations}\label{eqs:main}
\begin{equation}\label{1main}
\lim_{k\to \infty}\Big(\|u_k-u\|_{BD(\Omega)}+\|e(u_k) - e(u)\|_{L^p(\Omega;\Mnn)} + \hn(J_{u_k}\triangle J_u)\Big) = 0\,.
\end{equation}
Moreover, (if $p\in \big[1,\frac{n}{n-1}\big]$ this is trivial) there are Borel sets $E_k\subset \Omega$ such that
\begin{equation}\label{5main}
\lim_{k\to \infty}\mathcal{L}^n(E_k)= \lim_{k\to \infty}\int\limits_{\Omega\setminus E_k} |u_k -u|^p \dx = 0\,.
\end{equation}
\end{subequations}
\end{theorem}
The theorem above is sharp, in the sense that it provides the strongest possible approximation 
of all the relevant quantities in the definition of $SBD^p$. Moreover, differently from \cite{Cha04, Cha05Add} that assume $u\in L^2$, it does not require \emph{any additional integrability} assumption on $u$, and it is valid for any $p>1$ (in \cite{CFI17Density} it is observed that the construction in \cite{Cha04, Cha05Add} does not work for $p\neq 2$). These characteristics are in common with the sharp density result in $GSBD^p$ \cite{CC17}, which employs a similar construction, here improved to deal with $[u]$, see below.
 
 We remark that \cite{Iur14} and \cite{CC17} approximate also \emph{any truncation} of $[u]$, but this is not enough to deal with energies such as \eqref{eq:cohesive} without assuming \emph{a priori} a uniform $L^\infty$ bound.

It is interesting to compare Theorem~\ref{teo:density} with available density results in $SBV^p$, where of course there are more tools, such as the maximum principle or the coarea formula, due to the control on all $\nabla u$. On the one hand, Theorem~\ref{teo:density} may be combined with weaker $SBV^p$ approximations, but through functions with more regular jump set; on the other hand, our result provides stronger properties (some weaker) with respect to the available approximations in $BV$ norm for $SBV^p$, giving the possibility to improve them. 

First we consider the theorem by Cortesani and Toader, that approximates functions in $SBV^p\cap L^\infty$ with respect to an energy
\begin{equation}\tag{C'}\label{eq:enCorToa}
\int\limits_\Omega |\nabla u|^p \dx + \hn(J_u) + \int \limits_{J_u} \phi(x, u^+, u^-, \nu_u)  \dh\,,
\end{equation}
for very general $\phi$ (cf.\ Theorem~\ref{teo:CorToa}, see also the earlier \cite{DibSer97} for a weaker result, and \cite{AmarDeCicco} for an approximation for $BV\cap L^\infty$ functions). The approximating functions are of class $C^\infty \cap W^{m, \infty}$, for every $m\in \N$, outside the jump set, in turn closed and contained in a finite union of $(n{-}1)$-simplexes. Of course, this additional regularity on the jump set is in general in contrast to convergence in $BV$-norm. 

The first approximation result in $BV$-norm, for functions in $SBV^p \cap L^\infty$, is due to Braides and Chiadò-Piat \cite{BraChP96}: the approximating functions $u_k$ are $C^1$ outside some closed rectifiable sets $R_k$, such that $J_{u_k}\subset R_k$, with no information on the shape of $J_{u_k}$.

In the recent paper \cite{DPFusPra17}, De Philippis, Fusco, and Pratelli approximate $SBV^p$ functions 
by means of $u_k$ in $C^\infty(\Omega\sm J_{u_k})$, with $J_{u_k}$ 
 a compact $C^1$ manifold, up to a $\hn$-negligible set, and
\begin{equation*}
\lim_{k\to \infty}\Big(\|u_k-u\|_{BV(\Omega;\R^m)}+\|\nabla(u_k) - \nabla(u)\|_{L^p(\Omega;\mathbb{M}^{m{\times}n})} + \hn(J_{u_k}\triangle J_u)\Big) = 0\,.
\end{equation*}

The main improvement due to Theorem~\ref{teo:density}, besides the fact that it holds in $SBD^p$, is that our $u_k$ are also in $W^{m, \infty}(\Omega\sm J_{u_k})$, for every $m\in \N$, that may be very important in the applications. A possible weakness of our result is the fact that $J_{u_k}$ is not a $C^1$ manifold, even if, for the applications that we imagine at the moment (also for those presented in \cite{DPFusPra17}), 
one needs just $J_{u_k}$ closed, or one may employ \cite{CorToa99} (see also Remark~\ref{rem:1705181900}). 

In \cite{DPFusPra17} also two approximations in $BV$-norm, respectively for $SBV$ and $SBV^p_\infty$, are shown. In the spirit of this work, we prove the following approximations for $SBD$ and $SBD^p_\infty$. As in Theorem~\ref{teo:density}, we assume that $\Omega$ is open bounded Lipschitz. The crucial property is indeed that the trace of $u$ is integrable on $\partial\Omega$, so one could weaken the regularity assumption on $\Omega$.
\begin{theorem}\label{teo:TeorA}
Let $u\in SBD(\Omega)$. Then there exist $u_k\in SBD(\Omega)\cap L^\infty(\Omega;\Rn)$ such that $J_{u_k}$ is, up to a $\hn$-negligible set,
 a finite union of pairwise disjoint $C^1$ compact hypersurfaces  contained (strictly) in $\Omega$, $u_k\in C^\infty(\ol \Omega\sm J_{u_k};\Rn) \cap W^{m,\infty}(\Omega\sm J_{u_k};\Rn)$, and
\begin{equation}\label{eqs:mainTeorA}
\lim_{k\to \infty} \big(\|u_k-u\|_{BD(\Omega)}  + \hn(J_{u_k}\sm J_u) \big) = 0 \,.
\end{equation}
\end{theorem}

\begin{theorem}\label{teo:TeorB}
Let $u\in SBD^p_\infty(\Omega)$, with $p>1$. Then there exist $u_k\in SBV^p(\Omega;\Rn)\cap L^\infty(\Omega; \Rn)$ such that each
$J_{u_k}$ is closed 
and included in a finite union of closed connected pieces of $C^1$ hypersurfaces, $u_k\in C^\infty(\ol \Omega\sm J_{u_k};\Rn) \cap W^{m,\infty}(\Omega\setminus J_{u_k}; \Rn)$ for every $m\in \N$, and:
\begin{equation}\label{eqs:mainTeorB}
\lim_{k\to \infty}\Big(\|u_k - u\|_{BD(\Omega)} +\|e(u_k) - e(u)\|_{L^p(\Omega;\Mnn)} \Big) = 0\,.
\end{equation}
\end{theorem}

We observe that in Theorem~\ref{teo:TeorA} we have also the full regularity of $J_{u_k}$, so this in fact generalises \cite[Theorem~A]{DPFusPra17} allowing us to consider $SBD(\Omega)$ and $u_k$ of class $W^{m,\infty}$ outside $J_{u_k}$. (Indeed we employ \cite[Lemma~4.3]{DPFusPra17} to pass from $J_{u_k}$ included in, to $J_{u_k}$ essentially equal to the finite union of the desired $C^1$ hypersurfaces.) 

As for the approximation in $SBD^p_\infty$, we are not able to guarantee that $\hn(J_{u_k}\sm J_u)$ vanishes. This issue is also present in the corresponding \cite[Theorem~B]{DPFusPra17}, so Theorem~\ref{teo:TeorB} is not sharp (cf.\ Remark~\ref{rem:1805180933}).

In all the previous theorems, notice the strong convergence of $u_k$ to $u$ in $BD$ implies that (see \eqref{1205181701} and $|a||b|/\sqrt{2}\leq |a\odot b| \leq |a||b|$ for every $a$, $b$ in $\Rn$) 
\[
\hspace{-1em}\int \limits_{J_u \cup J_{u_k}} \hspace{-1em}\big|[u]-[u_k]\big| \dh \to 0\,.
\]

We conclude this introduction by briefly describing the proof strategy and possible applications of our results.

In all the three theorems, we assume $u$ extended with 0 outside $\Omega$, and we start from a set $\widehat{\Gamma} \in C^1$ with $\hn(\widehat{\Gamma} \sm J_u)$ and $\int_{J_u\sm \widehat{\Gamma}} \big| [u] \big| \dh$ small. In the spirit of \cite{Cha04}, we cover $\widehat{\Gamma}$ by cubes $Q_j$ splitted almost in two halves by this hypersurface, and we apply a rough approximation procedure in the complement of the union of the cubes, and in both sides of any cube with respect to $\widehat{\Gamma}$. 

We need different rough approximations for the $SBD^p$ and $SBD^p_\infty$ case, provided by Theorem~\ref{teo:rough} and Proposition~\ref{prop:roughTeorB}, respectively, while for Theorem~\ref{teo:TeorA} a suitable convolution is enough. The idea behind any rough approximation is to partition a given domain by cubes of sidelength $C k^{-1}$ and to detect the \emph{bad cubes}, i.e.\ those where the jump energy (or, similarly, the measure of the jump set for $SBD^p$) is not controlled well: in these cubes (indeed also in the adjacent \emph{boundary good cubes}) one sets $u_k$ as the infinitesimal rigid motion which is the ``mean'' of $u$, while in the remaining \emph{good cubes} one employs either a Korn-Poincaré-type inequality provided by \cite{CCF16} (cf.\ Proposition~\ref{prop:3CCF16}), or Lemma~\ref{le:controllo eu}, or a convolution with a radial kernel supported on a ball of radius $k^{-1}$, respectively. This construction differs from that of the rough approximation in \cite[Theorem~3.1]{CC17}, where $u_k= 0$ on the bad cubes, because we were there interested mainly to the measure of $J_{u_k}$, and not to control $[u_k]$.

A fundamental point is to separate the sets on which employ the rough approximation: first, this requires the function to be defined in a small neighbourhood of any subset,
to have the room for convolution; the second issue is to glue all the pieces obtained from the rough approximations in each subset. These problems could be solved by the technique in \cite{Cha04}, at the expense of assuming \emph{a priori} $u\in L^p$, since partitions of unity are needed, or by the trick in \cite{CC17}, that employs also an extension argument derived from Nitsche \cite{Nie81} (see Lemma~\ref{le:Nitsche}). Now there is a further delicate issue: if we glue as in \cite{CC17} we are not able to control $[u_k]$ on the intersection between $\partial Q_j$ and the zone where we extend by Lemma~\ref{le:Nitsche}, even if this has small $\hn$ measure. For this reason we have to perform a very careful approximation procedure, keeping the reflected zone of height $C k^{-1}$, so comparable to the size of small cubes and of the convolution kernels.

A key difference with respect to \cite{DPFusPra17} is that the rough approximants are smooth in a neighbourhood of any piece, so gluing them we keep the regularity up to the jump. This is not the case if one employs variable convolution kernels whose size decreases close to $\widehat{\Gamma}$, as in \cite{DPFusPra17}.

As application, we present an improvement to the sharp version of  two $\Gamma$-convergence approximations by \emph{phase-field} energies \emph{à la} Ambrosio-Tortorelli (cf.\ \cite{AmbTorCPAM}) for the energy \eqref{eq:cohesive}, in \cite{FocIur14} and \cite{CarVG18} (we mention also some approximations for cohesive energies \cite{CFI15}, \cite{DMOrlToa16}, \cite{BarLazZep16}). In \cite{FocIur14} and \cite{CarVG18}, the $\Gamma$-limsup inequality was proven just in $SBD^2 \cap L^\infty$, because this was done by hand for the regular functions provided by the Cortesani-Toader approximation, and then extended by \cite{Iur14}. Now it is enough to apply Theorem~\ref{teo:density} to pass directly to $SBD^2$, without any further integrability assumption.

We give no direct application to Theorems~\ref{teo:TeorA} and \ref{teo:TeorB}, but we recall that \cite[Theorem~6.1]{DPFusPra17} proves a representation formula for the total variation of $\mathrm{D}u$ for $BV$ and $SBV$ functions, derived from the analogous of Theorem~\ref{teo:TeorA} in \cite{DPFusPra17}.

In general, the result presented could be abstract tools useful to extend a variety of $\Gamma$-convergence approximations for e.g.\ suitable cohesive-type energies, that might be for instance in terms of finite elasticity or non-local energies, see respectively \cite{Fri17ARMA} and \cite{Neg03, Neg06} for the case of Griffith energy.

The plan of the paper is the following. In Section~\ref{Sec1} we fix the notation and recall some
technical lemmas, in Section~\ref{sec:densityrough} we present the rough approximation for Theorem~\ref{teo:density}, which is completely proven in Section~\ref{sec:proofMain}. Section~\ref{sec:other} is devoted to prove the other two density results, and the applications are contained in Section~\ref{sec:Appl}.

\section{Notation and preliminaries}\label{Sec1}

We denote by $\mathcal{L}^n$ and $\mathcal{H}^k$ the $n$-dimensional Lebesgue measure and the $k$-dimensional Hausdorff measure. For any locally compact subset $B$ of $\Rn$, the space of bounded $\R^m$-valued Radon measures on $B$ is indicated as $\mathcal{M}_b(B;\R^m)$. For $m=1$ we write $\mathcal{M}_b(B)$ for $\mathcal{M}_b(B;\R)$ and $\mathcal{M}^+_b(B)$ for the subspace of positive measures of $\mathcal{M}_b(B)$. For every $\mu \in \mathcal{M}_b(B;\R^m)$, $|\mu|(B)$ stands for its total variation.
We use the notation: $\chi_E$ for the indicator function of any $E\subset \R^n$, which is 1 on $E$ and 0 otherwise; $B_\varrho(x)$ for the open ball with center $x$ and radius $\varrho$; $x\cdot y$, $|x|$ for the scalar product and the norm in $\Rn$; $p^*$ for $np/(n-p)$, $n$ being the space dimension.
\par
\medskip
\paragraph{\bf $BV$ and $BD$ functions.}
For $U\subset \Rn$ open, a function $v\in L^1(U)$ is a \emph{function of bounded variation} on $U$, denoted by $v\in BV(U)$, if $\mathrm{D}_i v\in \mathcal{M}_b(U)$ for $i=1,\dots,n$, where $\mathrm{D}v=(\mathrm{D}_1 v,\dots, \mathrm{D}_n v)$ is its distributional gradient. A vector-valued function $v\colon U\to \R^m$ is $BV(U;\R^m)$ if $v_j\in BV(U)$ for every $j=1,\dots, m$.

%

The space of \emph{functions of bounded deformation} on $U$ is
\begin{equation*}
BD(U):=\{v\in L^1(U;\Rn) \colon \E v \in \mathcal{M}_b(U;\Mnn)\}\,,
\end{equation*}
where $\E v$ is the distributional symmetric gradient of $v$.
It is well known (see \cite{AmbCosDM97, Tem}) that for $v\in BD(U)$, the \emph{jump set} $J_v$, defined as the set of points $x\in U$ where $v$ has two different one sided Lebesgue limits $v^+(x)$ and $v^-(x)$ with respect to a suitable direction $\nu_v(x)\in \Sn$, is countably $(\hn, n-1)$ rectifiable (see, e.g.\ \cite[3.2.14]{Fed}), and that
\begin{equation*}
\mathrm{E}v=\mathrm{E}^a v+ \mathrm{E}^c v + \mathrm{E}^j v\,,
\end{equation*}
where $\mathrm{E}^a v$ is absolutely continuous with respect to $\mathcal{L}^n$, $\mathrm{E}^c v$ is singular with respect to $\mathcal{L}^n$ and such that $|\mathrm{E}^c v|(B)=0$ if $\hn(B)<\infty$, while 
\begin{equation}\label{1205181701}
\E^j v=[v]\odot \nu_v \,\hn  \mres  J_v\,.
\end{equation}
In the above expression of $\E^j v$, $[v]$ denotes the \emph{jump} of $v$ at any $x\in J_v$ and is defined by $[v](x):=(v^+-v^-)(x)$, the symbols $\odot$ and $\mres$ stands for the symmetric tensor product and the restriction of a measure to a set, respectively. Since $|a\odot b| \geq |a||b|/\sqrt{2}$ for every $a$, $b$ in $\Rn$, it holds $[v]\in L^1(J_v;\Rn)$. 
The density of $\mathrm{E}^a v$ with respect to $\mathcal{L}^n$ is denoted by $e(v)$, and we have that (see \cite[Theorem~4.3]{AmbCosDM97}) for $\mathcal{L}^n$-a.e.\ $x\in U$
\begin{equation*}\label{3105171931}
\lim_{\varrho\to 0^+}\frac{1}{\varrho^n} \int \limits_{B_\varrho(x)}\frac{\big(v(y)-v(x)-e(v)(x)(y-x)\big)\cdot (y-x)}{|y-x|^2}\mathrm{d}y =0\,.
\end{equation*}
The space $SBD(U)$ is the subspace of all functions $v\in BD(U)$ such that $\mathrm{E}^c v=0$, while for $p\in (1,\infty)$
\begin{equation*}
SBD^p(U):=\{v\in SBD(U)\colon e(v)\in L^p(U;\Mnn),\, \hn(J_v)<\infty\}\,.
\end{equation*}
Analogous properties hold for $BV$, as the countable rectifiability of the jump set and the decomposition of $\mathrm{D}v$.
Similarly, $SBV(U;\R^m)$ is the space of $BV(U;\R^m)$ with null Cantor part and
\begin{equation*}
SBV^p(U;\R^m):=\{v\in SBV(U;\R^m)\colon \nabla v \in L^p(U;\mathbb{M}^{m{\times}n}),\, \hn(J_v)<\infty\}\,,
\end{equation*}
 $\nabla v$ denoting the density of $\mathrm{D}^a v$, the absolutely continuous part of $\mathrm{D}v$, with respect to $\mathcal{L}^n$. Consider also the space (for this notation see e.g.\ \cite{DPFusPra17})
 \begin{equation*}
SBV^p_\infty(U;\R^m):=\{v\in SBV(U;\R^m)\colon \nabla v\in L^p(U;\mathbb{M}^{m{\times}n})\}\,,
\end{equation*}
and its analogous
 \begin{equation*}
SBD^p_\infty(U):=\{v\in SBD(U)\colon e(v)\in L^p(U;\Mnn)\}\,.
\end{equation*}
For more details on the spaces $BV$, $SBV$ and $BD$, $SBD$ we refer to \cite{AFP} and to \cite{AmbCosDM97, BelCosDM98, Bab15, Tem}, respectively.
Below we recall some other properties that will be useful in the following.

We start with an extension lemma derived from \cite[Lemma~1]{Nie81}. The result is employed in dimension 2 in \cite[Lemma~3.4]{CFI16ARMA}, and formulated in the more general setting of the space $GSBD^p$ in \cite[Lemma~5.2]{FriSol16} and in \cite[Lemma~2.8]{CC17}, to which we refer for more details of the proof.  
\begin{lemma}\label{le:Nitsche}
Let $R\subset \Rn$ be an open rectangle, $R'$ be the reflection of $R$ with respect to one face $F$ of $R$, and $\widehat{R}$ be the union of $R$, $R'$, and $F$. Let $v \in SBD^p(R)$. Then $v$ may be extended by a function $\widehat{v}\in SBD^p(\widehat{R})$ such that 
\begin{subequations}
\begin{align}
\hn(J_{\widehat{v}}&\cap F)=0\,,\label{2705170936}\\
\|\widehat{v}\|_{L^1(\widehat{R})} &\leq c \|v\|_{L^1(R)}\label{2004181331}\\
\hn(J_{\widehat{v}})&\leq c\, \hn(J_v)\,,\label{2705170937}\\
\int \limits_{J_{\widehat{v}}} |[\widehat{v}]|\dh &\leq c \int \limits_{J_v} |[v]|\dh\,,\label{2004180804}\\
\int\limits_{\widehat{R}} |e(\widehat{v})|^p\dx&\leq c\, \int\limits_R |e(v)|^p \dx\,,\label{2705170938}
\end{align}
\end{subequations}
for a suitable $c>0$ independent of $R$ and $v$.
\end{lemma}
\begin{proof}
We may follow \cite[Lemma~2.8]{CC17}, stated for $v\in GSBD^p(R)$. We assume that $F\subset \{(x',x_n)\in \R^{n-1}\times \R\colon x_n=0\}$  and $R \subset \{(x',x_n)\in \R^{n-1}\times \R\colon x_n<0\}$, fix any $\mu$, $\nu$ such that $0<\mu <\nu <1$, and let $q:=\frac{1+\nu}{\nu-\mu}$. 
Then $v'$ is defined on $R'$ by
\begin{equation*}
v':=q\, v_{A_\mu} + (1-q) v_{A_\nu}\,,
\end{equation*}
with $A_{\mu}=\mathrm{diag}\,(1,\dots,1,-\mu)$, $A_{\nu}=\mathrm{diag}\,(1,\dots,1,-\nu)$, and for any $u\in SBD^p(\Omega)$, $A \in M^{n{\times}n}$
\begin{equation}\label{2204182357}
u_A(x):=A^T u(Ax)\,.
\end{equation}
Following \cite[Lemma~2.8]{CC17}, it is immediate to verify that if $v\in SBD^p(R)$ then $\widehat{v}\in SBD^p(R')$ and \eqref{2705170936}, \eqref{2004181331}, \eqref{2705170937}, \eqref{2705170938} hold.
In order to show \eqref{2004180804} we notice that, for $u_A$ as in \eqref{2204182357}, $J_{u_A}=A^{-1}(J_u)$ and
\begin{equation*}
[u_A](A^{-1}x)=A^T[u](x)\,
\end{equation*}
for any $x\in J_u$.
This gives the further property corresponding to \cite[Lemma~2.7]{CC17} that allows us to repeat the argument of \cite[Lemma~2.8]{CC17} for the amplitude of the jump.
\end{proof}
We now recall the so called Korn-Poincaré inequality in $BD$ (cf.\ \cite{Koh82, Tem}). Notice that in the case of $W^{1,p}$ functions, with $p>1$, one obtains an analogous control for the $L^{p^*}$ norm of $u-a$ by combining the classical Korn and Poincaré inequalities.
\begin{proposition}\label{prop:KornPoinBD}
Let $U\subset \Rn$ be a bounded, connected, Lipschitz domain. Then there exists $c>0$ depending only on $U$ and invariant under rescaling of the domain, such that for every $u\in BD(U)$ there exists an affine function $a\colon \Rn\to\Rn$ with $e(a)=0$ such that
\begin{equation*}
\|u-a\|_{L^{1^*}(U;\Rn)} \leq c\, |\E u|(U)\,.
\end{equation*} 
In particular, for any cube $Q_r$ of sidelength $r$, H\"older inequality gives that
\begin{equation}\label{2704181803}
\|u-a\|_{L^1(Q_r;\Rn)}\leq c(Q_1)\, r \,|\E u|(Q_r)\,.
\end{equation}
\end{proposition}
Different Korn-Poincaré-type inequalities have been proven recently in the context of $SBD^p$. In \cite{CFI16ARMA, Fri17M3AS, FriPWKorn} also Korn-type inequalities have been considered. We recall here a result due to Chambolle, Conti, and Francfort and employed in \cite{CCF17, ChaConIur17, CC17, CC18}.
\begin{proposition}\label{prop:3CCF16}
Let $Q =(-r,r)^n$, $Q'=(-r/2, r/2)^n$, $u\in SBD^p(Q)$, $p\in [1,\infty)$, $\hn(J_u)<\infty$. Then there exist a Borel set $\omega\subset Q'$ and an affine function $a\colon \Rn\to\Rn$ with $e(a)=0$ such that $\mathcal{L}^n(\omega)\leq cr \hn(J_u)$ and
\begin{equation}\label{prop3iCCF16}
\int\limits_{Q'\setminus \omega}(|u-a|^{p}) ^{1^*} \dx\leq cr^{(p-1)1^*}\Bigg(\int\limits_Q|e(u)|^p\dx\Bigg)^{1^*}\,.
\end{equation}
If additionally $p>1$, then there is $q>0$ (depending on $p$ and $n$) such that, for a given mollifier $\varphi_r\in C_c^{\infty}(B_{r/4})\,, \varphi_r(x)=r^{-n}\varphi_1(x/r)$, the function $v=u \chi_{Q'\setminus \omega}+a\chi_\omega$ obeys
\begin{equation}\label{prop3iiCCF16}
\int\limits_{Q''}|e(v\ast \varphi_r)-e(u)\ast \varphi_r|^p\dx\leq c\left(\frac{\hn(J_u)}{r^{n-1}}\right)^q \int\limits_Q|e(u)|^p\dx\,,
\end{equation}
where $Q''=(-r/4,r/4)^n$.
The constant in (i) depends only on $p$ and $n$, the one in (ii) also on $\varphi_1$.  
\end{proposition}

 \begin{remark}\label{3005181123}
 By H\"older inequality and \eqref{prop3iCCF16} it follows that
 \begin{equation}\label{prop3CCFHolder}
\|u-a\|_{L^p(Q'\setminus \omega; \Rn)} \leq   c r \|e(u)\|_{L^p(Q;\Mnn)}\,.
 \end{equation}
 Moreover, looking at the proof of Proposition~\ref{prop:3CCF16} (take $g=|e(w)| \chi_Q$ instead of $g=|e(w)|^p \chi_Q$ and $p=1$ in the last part of \cite[Proposition~2]{CCF16}) one may see that for $a$ as in Proposition~\ref{prop:3CCF16} it holds also
  \begin{equation}\label{prop3CCFHolderesp1}
\|u-a\|_{L^1(Q'\setminus \omega; \Rn)} \leq   c r \|e(u)\|_{L^1(Q;\Mnn)}\,.
 \end{equation}
 \end{remark}
 
 In the following $\Omega$ will be a bounded open Lipschitz subset of $\Rn$.
We will denote by $C$ a generic positive constant depending only (at most) on $n$ and $p$, using $c$ only when we recall for the first time Lemma~\ref{le:Nitsche}, Proposition~\ref{prop:KornPoinBD} or Proposition~\ref{prop:3CCF16}.
\section{An auxiliary density result}\label{sec:densityrough}

\begin{theorem}\label{teo:rough}
Let $\Omega$, $\widetilde{\Omega}$ be bounded open subsets of $\Rn$, with $\ol \Omega\subset \widetilde{\Omega}$, $p\in (1,\infty)$, $\theta \in (0,1)$, and let $u\in SBD^p(\widetilde{\Omega})$.
Then there exist $u_k\in SBV^p(\Omega;\Rn)\cap L^\infty(\Omega; \Rn)$ such that $J_{u_k}$ is included in a finite union of $(n-1)$--dimensional closed cubes, $u_k\in C^\infty(\ol\Omega\setminus J_{u_k}; \Rn) \cap W^{m,\infty}(\Omega\setminus J_{u_k}; \Rn)$ for every $m\in \N$, and:
\begin{subequations}
\begin{align}
\limsup_{k\to \infty} \int\limits_\Omega  |e(u_k)|^p \dx & \leq \int\limits_\Omega |e(u)|^p \dx \,, \label{2rough}\\
\hn(J_{u_k}\cap \Omega)& \leq C \,\theta^{-1} \hn(J_u)\,,\label{3rough}\\
\limsup_{k\to \infty} \int\limits_{J_{u_k}} \big|[u_k]\big| \dh &\leq C \int\limits_{J_u} \big|[u]\big| \dh \,, \label{4rough}
\end{align}
for a suitable $C>0$ independent of $\theta$ and $k$.
Moreover, there are Borel sets $E_k\subset \Omega$ such that
\begin{equation}\label{5rough}
\lim_{k\to \infty}\mathcal{L}^n(E_k)= \lim_{k\to \infty}\int\limits_{\Omega\setminus E_k} |u_k -u|^p \dx = 0\,.
\end{equation}
In particular
\begin{align}
u_k &\wstar u \quad\text{in }BD(\Omega)\,,\label{1rough}\\
e(u_k)&\to e(u) \quad\text{in }L^p(\Omega;\Mnn)\,.\label{6rough}
\end{align}
\end{subequations}
\end{theorem}

\begin{proof}
As in \cite[Theorem~3.1]{CC17} we partition the domain into cubes of sidelength $k^{-1}$ and consider the cubes that contain a small amount of jump with respect to the perimeter of their boundary (in terms of the parameter $\theta$). While in these (good) cubes we do a construction as in \cite[Theorem~3.1]{CC17}, based on Proposition~\ref{prop:3CCF16},
we have to treat differently the remaining (bad) cubes. Indeed, even if we control in measure the perimeter of the union of the bad cubes, we have to define carefully the approximating functions in this zone in order to control the amplitude of the jump created on the perimeter.  
We then define in each bad cube with sidelength $k^{-1}$ (this is done also in cubes adjacent to bad cubes, called boundary good cubes) the $k$-th approximating function as the affine infinitesimal rigid motion given by Proposition~\ref{prop:KornPoinBD}:
in this way we introduce new jumps with respect to the construction in \cite[Theorem~3.1]{CC17}, but we estimate both their measure and the corresponding energy, in terms of the total variation of the symmetric gradient in (a neighbourhood of) the union of bad cubes.  As $k \to \infty$ one sees only the contribution of the jump part, since the $n$-dimensional measure of the union of bad cubes vanishes. 

We now recall notation and results from \cite[Theorem~3.1]{CC17}, and show the additional properties obtained by this different construction. In the following we omit to write the target spaces $\Rn$ or $\Mnn$ from the notation for the $L^p$ norm, to ease the reading.
\vspace{1em}

Let us fix an integer $k$ with $k> \frac{16 \sqrt{n}}{\mathrm{dist }(\partial \Omega, \partial \widetilde{\Omega})}$, let $\varphi$ be a smooth radial function with compact support in the unit ball $B(0,1)$, 
and let $\varphi_k(x)=k^n \varphi(kx)$. 
\newline
\\
{\bf Good and bad nodes.}
For any $z\in (2 \km) \Z^n \cap \Omega$ consider the cubes of center $z$
\begin{equation*}
\begin{split}
q_z^k&:=z+(-\km,\km)^n\,,\quad \hspace{1.3em}\tq_z^k:= z+(-2\km,2\km)^n\,,\\ Q_z^k&:=z+(-4\km,4\km)^n\,, \quad \widetilde{Q}_z^k:=z+(-8\km,8\km)^n\,.
\end{split}
\end{equation*}
The ``good'' and the ``bad'' nodes are defined as
\begin{equation}\label{eq:defGoodBad}
G^k:=\{z\in (2 \km) \Z^n \cap \Omega : \hn(J_u\cap \Qz)\leq \theta k^{-(n-1)}\}\,, \quad B^k:=(2 \km) \Z^n \cap \Omega \setminus G^k\,,
\end{equation}
to which correspond the subsets of $\widetilde{\Omega}$
\begin{equation}\label{1204181406}
\Omega_g^k:=\bigcup_{z\in G^k} \qz\,,\quad \widetilde{\Omega}^k_b:=\bigcup_{z\in B^k}\Qz\,.
\end{equation}
Notice that 
\begin{equation}\label{1205171038}
 \widetilde{\Omega}^k_b = \widetilde{\Omega}\setminus \Omega_g^k + (-3k^{-1},3k^{-1})^n\,,
\end{equation}
so that a row (and a half) of ``boundary'' cubes of $\Omega_g^k$ belongs to $\widetilde{\Omega}_b^k$ (see also the cubes in the second figure at page \pageref{fig}).
By \eqref{eq:defGoodBad}
\begin{equation}\label{1204181415}
\# B^k \leq \hn(J_u)\, k^{n-1} \theta^{-1}\,,
\end{equation}
and then
\begin{equation}\label{1805171025}
\mathcal{L}^n\left(\widetilde{\Omega}^k_b\right)\leq 16^n \hn(J_u)\, k^{-1}\,\theta^{-1}\,. 
\end{equation} 
Let us apply Proposition~\ref{prop:3CCF16} for any $z\in G^k$ (see also Remark~\ref{3005181123}).
Then there exist $\omega_z \subset \tqz$ and $a_z\colon \Rn \to \Rn$ affine with $e(a_z)=0$, such that (we recall directly only the condition corresponding to \eqref{prop3CCFHolder}, weaker than \eqref{prop3iCCF16}, and \eqref{prop3CCFHolderesp1})
\begin{equation}\label{1005171230}
\mathcal{L}^n(\omega_z)\leq c k^{-1} \hn(J_u\cap \Qz) \leq c \theta k^{-n}\,,
\end{equation}
\begin{equation}\label{prop3iCCF16applicata}
\|u-a_z\|_{L^{p}(\tqz\setminus \omega_z)} \leq ck^{-1} \|e(u)\|_{L^p(\Qz)}
\end{equation}
\begin{equation}\label{prop3iCCF16applicataesp1}
\|u-a_z\|_{L^{1}(\tqz\setminus \omega_z)} \leq ck^{-1} \|e(u)\|_{L^1(\Qz)}
\end{equation}
and
\begin{equation}\label{prop3iiCCF16applicata}
\begin{split}
\int\limits_{\qz}|e(v_z\ast \varphi_k)-e(u)\ast \varphi_k|^p\dx  \leq c\left(\hn(J_u\cap \Qz)\,k^{n-1}\right)^q \int\limits_{\Qz}|e(u)|^p\dx  \leq c \, \theta^q \int\limits_{\Qz}|e(u)|^p\dx \,,
\end{split}
\end{equation}
for $v_z:= u\chi_{\tqz\setminus \omega_z}+a_z \chi_{\omega_z}$ and a suitable $q>0$ depending on $p$ and $n$.

We define
\begin{equation*}
\omega^k:= \bigcup_{z\in G^k} \omega_{z}\,,\qquad
E_k:=\widetilde{\Omega}_b^k \cup \omega^k\,.
\end{equation*}
By \eqref{1005171230} we have 
\[
\mathcal{L}^n(\omega^k)\leq c k^{-1} \sum_{z\in G^k}  \hn(J_u\cap Q_{z}^k) \leq c \hn(J_u)\, k^{-1}\,,
\] so that \eqref{1805171025} implies
\begin{equation}\label{1805171033}
\lim_{k\to \infty}\mathcal{L}^n(E_k)=0\,.
\end{equation}

For every $z \in (2k^{-1})\mathbb{Z}^n \cap \Omega$ we employ Proposition~\ref{prop:KornPoinBD} and let $\tilde{a}_z\colon \Rn \to \Rn$ be the affine function with $e(\tilde{a}_z)=0$ such that (also here we recall directly \eqref{2704181803})
\begin{equation}\label{eq:PoinBad}
\| u-\tilde{a}_z \|_{L^1(\tilde{q}_z)} \leq C k^{-1} | \E u|(\tqz)\,.
\end{equation}
We remark that for every $z\in G^k$
\begin{equation}\label{1204180938}
\mathcal{L}^n(\tilde{q}_z)\|a_z-\tilde{a}_z\|_{L^\infty(\tqz)}\leq Ck^{-1} \Big(|\E u|(\tilde{q}^k_z)+\|e(u)\|_{L^1(\Qz)}\Big) \,.
\end{equation}
Indeed, by \eqref{prop3iCCF16applicataesp1} and
\eqref{eq:PoinBad} we get
\begin{equation}\label{1204180954}
\|a_z-\tilde{a}_z\|_{L^1(\tqz\setminus \omega_z)}  \leq Ck^{-1} \Big(|\E u|(\tqz)+\|e(u)\|_{L^1(\Qz)}\Big)\,,
\end{equation}
and then we deduce \eqref{1204180938} because
\begin{equation*}
\mathcal{L}^n(\tqz)\|a_z-\tilde{a}_z\|_{L^\infty(\tqz)} \leq C \|a_z-\tilde{a}_z\|_{L^1(\tqz\setminus \omega_z)}\,,
\end{equation*}
which is obtained following the argument of \cite[Lemma~4.3]{ConFocIur15} (see also \cite[Lemma~2.12]{CC17}), since $a_z-\tilde{a}_z$ is affine and $\mathcal{L}^n(\omega_z) \leq \mathcal{L}^n(\tqz)/4$.

%
%

\vspace{1em}
{\bf The approximating functions.}
Let $G^k=(z_j)_{j\in J}$, so that we order (arbitrarily) the elements of $G^k$,
 and define  
\begin{equation}\label{eq:defappr1}
\widetilde{u}_k:=
\begin{cases}
 u \quad &\text{in }\widetilde{\Omega}\setminus \omega^k\,,\\
 a_{z_j}\quad &\text{in }\omega_{z_j}\setminus \bigcup_{i<j}\omega_{z_i}\,,
\end{cases}
\end{equation}
and
\begin{equation}\label{eq:defapprox}
u_k:= \begin{cases}
\widetilde{u}_k \ast \varphi_k \quad &\text{in }\Omega\setminus \widetilde{\Omega}_b^k\,,\\
\tilde{a}_z \quad &\text{in } \qz \cap \widetilde{\Omega}_b^k\,.
\end{cases}
\end{equation}
It is immediate that $u_k \in SBV^p(\Omega;\Rn)\cap L^\infty(\Omega;\Rn)$, since $u \in BD(\Omega)\subset L^1(\Omega;\Rn)$, and that $u_k\in C^\infty(\ol\Omega\sm J_{u_k};\Rn) \cap W^{m,\infty}(\Omega\sm J_{u_k};\Rn)$ for every $m\in \N$, since $\widetilde{u}_k\ast \varphi_k$ is smooth in a neighbourhood of $\Omega\setminus \widetilde{\Omega}_b^k$. Moreover $J_{u_k}$ is closed and included in a finite union of boundaries of $n$-dimensional cubes $\qz$.
\newline
\\
{\bf Proof of \eqref{3rough}.}
We have that 
\begin{equation*}
J_{u_k}\subset \overline{\widetilde{\Omega}_b^k}\,,
\end{equation*}
so the definition \eqref{1204181406} of $\widetilde{\Omega}_b^k$ gives
\begin{equation}\label{1204181419}
J_{u_k}\subset \bigcup_{z\in B^k} (J_{u_k}\cap \ol\Qz)\,.
\end{equation}
Notice that for every $\widehat{z}\in B^k$
\[
J_{u_k}\cap \ol Q^k_{\widehat{z}}= \partial Q^k_{\widehat{z}} \cup \bigcup_{\qz \subset Q^k_{\widehat{z}}} \partial \qz\,,
\]
and then
\begin{equation}\label{1204181414}
\hn(J_{u_k}\cap \ol Q^k_{\widehat{z}}) \leq C k^{-(n-1)}\,.
\end{equation}
for $C$ depending only on $n$.
Together with \eqref{1204181415} and \eqref{1204181419}, \eqref{1204181414} implies \eqref{3rough}.
\newline
\\
{\bf Proof of \eqref{4rough}.} In order to prove \eqref{4rough} we estimate the amplitude of the jump in two different sets: the common boundaries between cubes of sidelength $2k^{-1}$ included in $\widetilde{\Omega}_b^k$ (which give the jump of $u_k$ included in the interior of $\widetilde{\Omega}_b^k$) and $\partial \widetilde{\Omega}_b^k$, which is essentially (up to a $\hn$-negligible set) contained in the interior of suitable cubes of sidelength $2k^{-1}$, recall \eqref{1205171038}.

Let $q^k_{z}$ and $q^k_{z'}$ be included in $\widetilde{\Omega}_b^k$, with $\hn(\partial q^k_{z} \cap \partial q^k_{z'})>0$. Then \eqref{eq:PoinBad} gives
\begin{equation}\label{1204182142}
\begin{split}
\|\tilde{a}_z-\tilde{a}_{z'}\|_{L^1(\tqz \cap \tilde{q}^k_{z'})} \leq  \|u-\tilde{a}_z\|_{L^1(\tqz \cap \tilde{q}^k_{z'})}+\|u-\tilde{a}_{z'}\|_{L^1(\tqz \cap \tilde{q}^k_{z'})} \leq C k^{-1}|\E u|(\tilde{q}^k_z\cup \tilde{q}^k_{z'}) \,.
\end{split}
\end{equation}
Being $\tilde{a}_z-\tilde{a}_{z'}$ affine, we have that
\begin{equation*}
\frac{4^n}{2}k^{-n}\|\tilde{a}_z-\tilde{a}_{z'}\|_{L^\infty(\tqz \cap \tilde{q}^k_{z'})}= \mathcal{L}^n(\tqz \cap \tilde{q}^k_{z'}) \|\tilde{a}_z-\tilde{a}_{z'}\|_{L^\infty(\tqz \cap \tilde{q}^k_{z'})}\leq C \|\tilde{a}_z-\tilde{a}_{z'}\|_{L^1(\tqz \cap \tilde{q}^k_{z'})} \,,
\end{equation*}
and together with \eqref{1204182142} this gives
\begin{equation*}
\begin{split}
\int \limits_{\partial q^k_{z} \cap \partial q^k_{z'}} \hspace{-1em}\big|[u_k]\big| \dh &= \hspace{-1em}\int \limits_{\partial q^k_{z} \cap \partial q^k_{z'}} \hspace{-1em}|\tilde{a}_z-\tilde{a}_{z'}| \dh \leq 2^{n-1}k^{-(n-1)} \|\tilde{a}_z-\tilde{a}_{z'}\|_{L^\infty(\tqz \cap \tilde{q}^k_{z'})}\\
& \leq C\, |\E u|(\tilde{q}^k_z\cup \tilde{q}^k_{z'}) \,. 
\end{split}
\end{equation*}
We put together all these contributions, observing that the cubes $\tqz$ are finitely overlapping and $\tqz \subset \widetilde{\Omega}_b^k$  if $q_z\subset \widetilde{\Omega}_b^k$ (cf.\ \eqref{1205171038}). We therefore obtain that
\begin{equation}\label{1304180952}
\int \limits_{\interior{(\widetilde{\Omega}_b^k)}} \big| [u_k] \big| \dh \leq C\,|\E u|(\widetilde{\Omega}_b^k)\,.
\end{equation}

Let us now consider a node $z$ such that $q_z \cap \partial \widetilde{\Omega}_b^k \neq \emptyset$. By definition of $\widetilde{\Omega}_b^k$ we have that $z \in G^k\cap  \partial \widetilde{\Omega}_b^k$. We claim that
\begin{equation}\label{1204182302}
\|\widetilde{u}_k-a_z\|_{L^1(\tqz)} \leq C k^{-1}  \|e(u)\|_{L^1(\widetilde{Q}^k_z)} \,.
\end{equation}
Indeed \eqref{prop3iCCF16applicataesp1}
and the fact that $\omega_z \subset \omega^k$ implies that (recall that $u_k=u$ in $\tqz\sm \omega^k$ by definition)
\begin{equation*}
\|\widetilde{u}_k-a_z\|_{L^1(\tqz\sm \omega^k)} \leq C k^{-1}  \|e(u)\|_{L^1(\Qz)}\,,
\end{equation*}
and it is proven in \cite[equation (3.19)]{CC17} (the definition of $\widetilde{u}_k$ is the same, take in \cite[equation (3.19)]{CC17} the 
version with $p=1$) that
\begin{equation*}\label{1204182325}
\|\widetilde{u}_k-a_z \|_{L^1(\tqz\cap \omega^k)} \leq C\, \theta k^{-1} \|e(u)\|_{L^1(\widetilde{Q}^k_z)}\,,
\end{equation*}
thus \eqref{1204182302} is proven.

We now combine \eqref{1204182302} with \eqref{1204180938}, giving
\begin{equation*}
\|a_z-\tilde{a}_z\|_{L^1(\tqz)} \leq Ck^{-1} \Big(|\E u|(\tilde{q}^k_z)+ \|e(u)\|_{L^1(\Qz)}\Big) \,,
\end{equation*}
to get
\begin{equation}\label{1204182345}
\|\widetilde{u}_k-\tilde{a}_z\|_{L^1(\tqz)} \leq Ck^{-1} \Big(|\E u|(\tilde{q}^k_z)+\|e(u)\|_{L^1(\widetilde{Q}^k_z)}\Big) \,.
\end{equation}

It follows that for every $x \in \partial \widetilde{\Omega}_b^k \cap \qz$
\begin{equation*}
\begin{split}
\big|[u_k]\big|(x)=|u_k-\tilde{a}_z|(x) &\leq C \|\varphi\|_{L^\infty(B_1)} k^{n} \| \widetilde{u}_k-\tilde{a}_z\|_{L^1(B_{k^{-1}}(x))}\\&\leq C k^{n} \| \widetilde{u}_k-\tilde{a}_z\|_{L^1(\tqz)}\leq Ck^{n-1} \Big(|\E u|(\tilde{q}^k_z)+\|e(u)\|_{L^1(\widetilde{Q}^k_z)}\Big) \,,
\end{split}
\end{equation*}
where we used the fact that $\varphi_k\ast \tilde{a}_z=\tilde{a}_z$, being $\varphi$ radial and $\tilde{a}_z$ affine. We then conclude
\begin{equation}\label{1304180924}
\int \limits_{\partial \widetilde{\Omega}_b^k \cap \qz} \big|[u_k]\big| \dh \leq C \Big(|\E u|(\tilde{q}^k_z)+\|e(u)\|_{L^1(\widetilde{Q}^k_z)}\Big)\,.
\end{equation}
Let us sum up over $z\in G^k$ such that $\hn(\partial \widetilde{\Omega}_b^k \cap \qz)>0$, namely over $z\in G^k \cap \partial \widetilde{\Omega}_b^k$. We remark that
%
\begin{equation*}
\bigcup_{z\in G^k \cap  \partial \widetilde{\Omega}_b^k} \hspace{-1em}\tqz \subset \bigcup_{z' \in B^k} z'+(-6k^{-1}, 6 k^{-1})^n\,,\quad \bigcup_{z\in G^k \cap  \partial \widetilde{\Omega}_b^k} \hspace{-1em}\widetilde{Q}^k_z \subset \bigcup_{z' \in B^k} z'+(-12 k^{-1}, 12 k^{-1})^n=: \widetilde{\Omega}^k_{b,1}\,.
\end{equation*}
 Moreover the cubes $\tqz$, $\widetilde{Q}^k_z$ are finitely overlapping, and then by \eqref{1304180924} we deduce that
\begin{equation}\label{1304180953}
\int \limits_{\partial \widetilde{\Omega}_b^k} \big| [u_k] \big| \dh \leq C \int \limits_{J_u} \big| [u] \big| \dh + C \int \limits_{\widetilde{\Omega}^k_{b,1}} |e(u)| \dx\,.
\end{equation}
Collecting \eqref{1304180952} and \eqref{1304180953} we get (recall the definition of $u_k$) 
\begin{equation*}
\int \limits_{J_{u_k}} \big| [u_k] \big| \dh \leq C \int \limits_{J_u} \big| [u] \big| \dh + C \int \limits_{\widetilde{\Omega}^k_{b,1}} |e(u)| \dx\,.
\end{equation*}
By \eqref{1204181415} we get, as in \eqref{1805171025}, that $\mathcal{L}^n(\widetilde{\Omega}^k_{b,1}) \leq C \hn(J_u) k^{-1} \theta^{-1}$, so we conclude \eqref{4rough}.
\newline
\\
{\bf Proof of the remaining properties.} We notice that our definition of $u_k$ differs form that one in \cite[Theorem~3.1]{CC17} only in $\widetilde{\Omega}_b^k$, since there the approximating functions were set equal to 0. In particular we may employ properties referring to cubes in $\Omega\sm \widetilde{\Omega}_b^k$ proven in \cite[Theorem~3.1]{CC17}.

Combining \cite[equations (3.14), (3.15), (3.16), (3.19)]{CC17} we have that
\begin{equation}\label{2004181946}
\|u-u_k\|_{L^p((\Omega\sm \widetilde{\Omega}_b^k)\sm \omega^k)}\leq C k^{-1}\|e(u)\|_{L^p(\widetilde{\Omega})}\,.
\end{equation}
Moreover we may follow the argument to prove property (3.1d) in \cite{CC17} (with $\psi=|\cdot|$, now $(\mathrm{HP}\psi)$ are useless) to get
\begin{equation}\label{2004181910}
\|u-u_k\|_{L^1((\Omega\sm \widetilde{\Omega}_b^k)\cap \omega^k)}\leq C k^{-1}\|e(u)\|_{L^1(\widetilde{\Omega})}+C \theta \|u\|_{L^1(\widetilde{\Omega}_{g,2}^k)} + C k^{-1/2} \|u\|_{L^1(\widetilde{\Omega})}\,.
\end{equation}
The set $\widetilde{\Omega}_{g,2}^k$ above is defined as follows: 
we set $G^k_1$ as the good nodes for which the condition on $J_u$ is satisfied for $k^{-\frac{1}{2}}$ in place of $\theta$
\begin{equation*}
G^k_1:=\{z\in G^k \colon \hn(J_u\cap \Qz)\leq k^{-(n-\frac{1}{2})}\}\,, \qquad G^k_2:= G^k\setminus G^k_1\,.
\end{equation*}
and the set $\widetilde{G}^k_1$ of the nodes adjacent to nodes in $G^k_1$
\begin{equation*}
\begin{split}
\widetilde{G}^k_1&:=\{z\in G^k \colon \ol z \in G^k_1 \text{ for each } \ol z \in (2k^{-1}) \Z^n \text{ with } \| z - \ol z\|_\infty= 2k^{-1}\}\,,\\
\widetilde{G}^k_2&:=\{z\in G^k \colon \text{ there exists } \ol z \in G^k_2 \text{ with } \| z - \ol z\|_\infty= 2k^{-1}\}\,,
\end{split}
\end{equation*}
and then
\begin{equation*}\label{2004181930}
\widetilde{\Omega}^k_{g,2}:= \bigcup_{z_j \in \widetilde{G}^k_2} \widetilde{Q}_{z_j}\,.
\end{equation*}
We get that $\# G^2_k \leq \hn(J_u)\, k^{n-\frac{1}{2}}$, so 
\begin{equation*}
\# \widetilde{G}^k_2\leq (3^n-1) \hn(J_u)\, k^{n-\frac{1}{2}}\,.
\end{equation*}
In particular,
\begin{equation}\label{2305181808}
\mathcal{L}^n(\widetilde{\Omega}^k_{g,2})\leq C\, \hn(J_u) k^{-\frac{1}{2}}\,.
\end{equation}

Furthermore, the definition \eqref{eq:defapprox} of $u_k$ and \eqref{eq:PoinBad} give
\begin{equation}\label{2004181940}
\|u-u_k\|_{L^1(\widetilde{\Omega}^k_b)}\leq C k^{-1} |\E u|\big(\widetilde{\Omega}^k_b + (-2 k^{-1}, 2 k^{-1})^n\big)\,.
\end{equation}

Since it is still true that $e(u_k)=0$ on $\widetilde{\Omega}_b^k$ because $e(\tilde{a}_z)=0$, we get for free \eqref{2rough}, that is \cite[property (3.1b)]{CC17}. More precisely, by \cite[eqs. (3.32), (3.33)]{CC17} we have
\begin{equation}\label{2004182216}
\|e(u_k)\|_{L^p(\Omega)}\leq (1+C k^{-q})\|e(u)\|_{L^p(\widetilde{\Omega})}+C\, \theta^q \,\|e(u)\|_{L^p(\widetilde{\Omega}^k_{g,2})}\,,
\end{equation}
for $q>0$ depending only on $p$ and $n$.

Collecting \eqref{1805171033}, \eqref{2004181946}, \eqref{2004181910}, and \eqref{2004181940} we obtain \eqref{5rough} 
and $u_k \to u$ in $L^1(\Omega)$.
We have proven in particular that $u_k$ is bounded in $BD(\Omega)$, so the $L^1$ convergence to $u$ implies \eqref{1rough}, and \eqref{6rough} follows immediately from \eqref{2rough} (recall \cite[Theorem~1.1]{BelCosDM98}).
This concludes the proof.
\end{proof}

\section{Proof of the main density theorem}\label{sec:proofMain}

\begin{proof}[Proof of Theorem~\ref{teo:density}]

As in \cite[Theorem~1.1]{CC17}, the starting point is to cover most of $J_u$ and of $\partial\Omega$ by cubes for which $J_u$ or $\partial\Omega$ is almost a diameter, namely these sets
are there close (with respect to $\hn$ measure) to an almost flat $C^1$ hypersurface. The idea, introduced first in \cite{Cha04}, is to apply then the rough approximation on the one hand in both the (almost) half cubes in which the flat hypersurface splits each cube, and on the other hand in the remaining part of $\Omega$, since in all these sets the amount of jump is small.

We now recall the covering obtained in the first part of \cite[Theorem~1.1]{CC17}, referring to that theorem for details.
\newline
\\
{\bf Approximation of $J_u$ and $\partial \Omega$.}
For every $\varepsilon>0$, there exist a finite family of pairwise disjoint closed cubes $(\ol{Q_j})_{j=1}^{\ol{\jmath}}\subset \Omega$ with \[\ol{Q_j}=\ol Q(x_j,\varrho_j) \quad \text{for } x_j \in J_u\,\text{ and one face of $\ol{Q_j}$ normal to $\nu_u(x_j)$}\,,\]
$\nu_u(x_j)$ denoting the normal to $J_u$ at $x_j$,
and $C^1$ hypersurfaces $(\Gamma_j)_{j=1}^{\ol{\jmath}}$ with $x_j\in\Gamma_j$
such that
 \begin{subequations}\label{eqs:2405171202}
\begin{align} 
&\hn\Big(J_u \setminus \bigcup_{j=1}^{\overline{\jmath}}Q_j\Big)<\varepsilon\,, \label{1305171147}\\
 \hn\big((J_u\triangle \Gamma_j)\,\cap &\, \ol{Q_j}\big)< \varepsilon (2\varrho_j)^{n-1}< \, \frac{\varepsilon}{1-\varepsilon}  \hn(J_u\cap \ol{Q_j})\,, \label{1305171150}\\
 \Gamma_j\, \text{is a $C^1$ graph with respect} & \text{ to }  \nu_u(x_j) \text{ with Lipschitz constant less than $\varepsilon/2$}\,.\label{1904181025}
 \end{align}
 \end{subequations}
 In particular, \eqref{1904181025} gives
 \[\Gamma_j \subset \Big\{ x_j+ \sum_{i=1}^{n-1}  y_i\,b_{j,i} + y_n\, \nu_u(x_j) \colon  y_i \in (-\varrho_j,  \varrho_j),\,  y_n \in \Big(-\frac{\varepsilon \varrho_j}{2}, + \frac{\varepsilon \varrho_j}{2}\Big)\Big\} \,, 
 \]
  where 
 $
 (b_{j,i})_{i=1}^{n-1}\quad\text{ is an orthonormal basis of }\nu_u(x_j)^\perp$.
 
Arguing similarly for $\partial\Omega$ in place of $J_u$, there exist a finite family of closed cubes $(\ol Q_h^0)_{h=1}^{\ol h}$ of centers $x_h^0\in \partial\Omega$ and sidelength $2\varrho_h^0$, with one face normal to $\nu_\Omega(x_h^0)$ (the outer normal to $\Omega$ at $x_h^0$), pairwise disjoint and with empty intersection with any $\ol Q_j$, and $C^1$ hypersurfaces $(\Gamma_h^0)_{h=1}^{\ol h}$ with $x_h^0\in \Gamma_h^0$, such that
\begin{subequations}\label{eqs:0707172255}
\begin{align}
\hn\Big(&\partial  \Omega\setminus  \bigcup_{h=1}^{\ol h} Q_h^0\Big) <\varepsilon\,,\label{2105171956}\\
\hn\big((\partial\Omega\triangle \Gamma_h^0)\cap \ol{Q}_h^0 \big)&<\varepsilon (2\varrho_h^0)^{n-1}  < \frac{\varepsilon}{1-\varepsilon}\hn(\partial\Omega\cap \ol{Q}_h^0)\,,\label{2505171242}\\
\Gamma_h^0\, \text{is a $C^1$ graph with respect}  \text{ to }  &\nu_\Omega(x_h^0) \text{ with Lipschitz constant less than $\varepsilon/2$}\,.
\end{align}
\end{subequations}
 Notice that
 we may assume that conditions \eqref{eqs:2405171202} and \eqref{eqs:0707172255} hold also for 
 the enlarged cubes
 \[
 \ol Q_j + (-t,t)^n\,,\qquad \ol Q_h^0 + (-t,t)^n\,,
 \]
 for $t$ much smaller than $\varepsilon$ and $\min_{j,h}\{\varrho_j, \varrho_h^0\}$ (we will consider below a parameter $k$ chosen such that $k^{-1}$ is much smaller than $t$).
 
 We denote
\begin{equation}\label{2204181215}
\widehat{\Gamma}:= \bigcup_{j=1}^{\ol \jmath} (Q_j\cap \Gamma_j)\,, \qquad \qquad\widehat{\Gamma}_{\partial\Omega}:=\bigcup_{h=1}^{\ol h} (Q_h^0\cap \Gamma_h^0)\,. 
\end{equation} 
From \eqref{1305171147}, \eqref{1305171150}, and \eqref{2105171956}, \eqref{2505171242} it follows that 
\begin{equation}\label{2204181223}
\hn(J_u \triangle \widehat{\Gamma}) < C\, \varepsilon\, \hn(J_u)\,,\qquad \hn(\partial\Omega \triangle \widehat{\Gamma}_{\partial\Omega}) < C\, \varepsilon\, \hn(\partial\Omega)\,.
\end{equation}
 Let
\begin{equation}\label{2204181256}
\eta_\varepsilon:= \varepsilon \vee \Big(\int \limits_{J_u \sm \widehat{\Gamma}} \big|[u]\big| \dh + \int\limits_{\partial\Omega\sm \widehat{\Gamma}_{\partial\Omega}} |\mathrm{tr}_\Omega u|\dh \Big)^{1/(n-1)} \,.
\end{equation}
Then $\lim_{\varepsilon\to 0} \eta_{\varepsilon} =0$, since $[u] \in L^1(J_u; \Rn)$ and $\mathrm{tr}_\Omega u\in L^1(\partial\Omega;\Rn)$, being $\Omega$ Lipschitz and $u\in SBD(\Omega)$.
Moreover, we set
\begin{equation}\label{eq:defb0}
B_0:= \Omega\setminus \Big(\bigcup_{j=1}^{\overline{\jmath}} \ol{Q_j} \cup \bigcup_{h=0}^{\ol h} \ol Q_h^0\Big)\,.
\end{equation}
\newline
\\
{\bf Definition of the approximating functions in the cubes.}
We now describe the construction of the approximating functions in each cube of the collection $(\ol{Q_j})_{j=1}^{\ol{\jmath}}$ or $(\ol Q_h^0)_{h=1}^{\ol h}$. We then fix a single cube that we denote by $\ol Q= \ol Q(x, \varrho)$ and we call $\Gamma$ the corresponding hypersurface that splits $Q$ in two (almost) half cubes
$Q^+$ and $Q^-$, to ease the reading ($\Gamma$ is either close to $J_u$ or to $\partial\Omega$). We also assume that $x=0$ and $\nu(x)=e_n$. 

As in the case of the rough approximation in Theorem~\ref{teo:rough}, also now a construction finer than the corresponding one in \cite{CC17} is needed. In \cite{CC17} one constructs an auxiliary function in a neighbourhood of both the half cubes in a single step, employing
a unique extension for each half cube: in the strip of height $\varepsilon\varrho$ containing the jump the original function $u$ was replaced employing values of $u$ in the strip of the same size which is immediately below (for $Q^-$) or above (for $Q^+$). The argument in \cite[Theorem~1.1]{CC17} continues by applying the rough approximation to the auxiliary functions in both the half cubes and in $B_0$ and gluing simply by characteristic functions. In this way one introduces a further jump: even if its $\hn$-measure is $C\varepsilon \varrho^{n-1}$ at the boundary of each cube (so the total surface is small), its amplitude is unfortunately not controlled.

The idea is now to modify the original function on a strip of height $k^{-1}$ around $\Gamma$ in order to construct the approximation $u_k$ in each half cube, since this works well with convolution with kernels supported on $B(0,k^{-1})$. One has to choose carefully the zone where the function is extended from the two sides of $\Gamma$, in order to control the $\hn$-measure of the new jump set.


%


\begin{figure}[h]\label{fig}
\hspace{-2em}
\begin{minipage}[c]{0.49\linewidth}
\includegraphics[width=\linewidth]{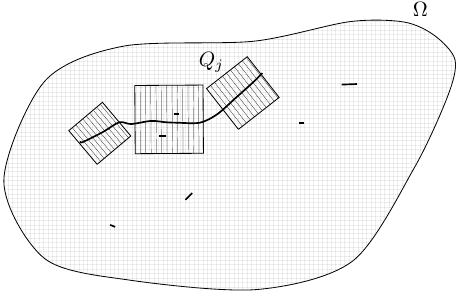}
\end{minipage}
\hfill
\begin{minipage}[c]{0.53\linewidth}
\includegraphics[width=\linewidth]{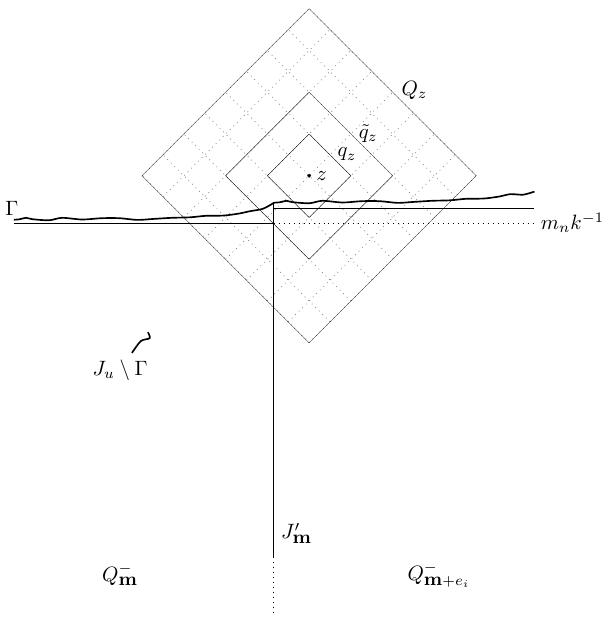}
\end{minipage}
\caption*{In the first figure, the cubes $Q_j$ covering almost all $J_u$, the relative (almost) parallelepipeds $Q^\pm_\m$, and the small cubes $\qz$ partitioning $\Omega$. In the second one, a zoom on the zone between two adjacent (almost) parallelepipeds for a cube $Q_j$, oriented with $\nu_u(x_j)$: $J_{u_k}$ is included in the union of the $\Gamma_j$, of the $J'_\m$ and of the boundary of bad cubes and boundary good cubes of sidelength $2 k^{-1}$.}
\end{figure}

In order to perform the following construction also for the cubes $(Q_h^0)^\pm$ (and after for $B_0$), that are possibly not included in $\Omega$, we extend 
$u$ outside $\Omega$ with the value 0.

Let $k\in \N$ be much larger than $(\eta_\varepsilon \varrho)^{-1}$, let us fix the (almost) half cube $Q^-$ and partition $Q^-$ into the union of (almost) parallelepipeds. We denote 
\begin{align*}
& F_{\textbf{m}}:=\big\{ (y_1,\dots, y_{n-1})\in \R^{n-1}\colon y_i \in (\eta_\varepsilon k)^{-1} m_i + \big(0, (\eta_\varepsilon k)^{-1})\big\}\quad\\
&F'_{\textbf{m}}:= F_{\textbf{m}} + (-32\sqrt{n} k^{-1},32 \sqrt{n} k^{-1})^{n-1}  \,,
\end{align*}
where $F$ stands for ``face'' of an $n$-dimensional cube, for (we may assume $\eta_\varepsilon k \varrho \in \N$, otherwise below put the fractional part of $-\eta_\varepsilon k \varrho$ in place of 0)
\begin{equation}\label{2904181212}
\textbf{m}=(m_1,\dots,m_{n-1})\in \{-\eta_\varepsilon k \varrho, -\eta_\varepsilon k \varrho+1, \dots, 0, \dots, \eta_\varepsilon k \varrho-1\}^{n-1} \subset \N^{n-1}\,.
\end{equation}
Since $\Gamma$ is the graph of a $\varepsilon/2$-Lipschitz function with respect to $e_n$ and $\eta_\varepsilon\geq \varepsilon$, there exists $m_n \in \R$, depending on $\textbf{m}$, such that
\begin{equation}\label{1904182312}
\Gamma\cap \big(F'_{\textbf{m}} {\times} (-\varrho, (\varepsilon \varrho)/2)\big) \subset F'_{\textbf{m}} {\times} (m_n, m_n+1/2)k^{-1}\,, 
\end{equation}
where $ (m_n, m_n+1/2)k^{-1}=  (m_n k^{-1}, (m_n+1/2)k^{-1}) \subset \R$ (indeed every side of $F_\m$ has length $\eta_\varepsilon^{-1}k^{-1}\leq \varepsilon^{-1}k^{-1}$).
Let us set
\begin{equation}\label{2004180729}
u_{\textbf{m}}:=\begin{dcases}
u \quad &\text{in } F'_{\m}{\times}(-\varrho-16\sqrt{n} k^{-1}, m_n k^{-1})\\
\widehat{u} \quad &\text{in }F'_{\m}{\times}\big((m_n,m_n+25\sqrt{n})k^{-1}\big)\,,
\end{dcases}
\end{equation}
where $\widehat{u}$ is obtained by Lemma~\ref{le:Nitsche} taking $F'_{\m}{\times}\{m_n k^{-1}\}$, 
\begin{equation}\label{2104181322}
R_{\m}:=F'_{\m}{\times}\big((m_n-25\sqrt{n},m_n)k^{-1}\big)\,,\qquad R'_{\m}:=F'_{\m}{\times}\big((m_n,m_n+25\sqrt{n})k^{-1}\big)
\end{equation}
 as $F$, $R$, $R'$ therein, respectively.
We introduce (see figures at page \pageref{fig})
\begin{equation*}
Q^-_{\m}:=Q^-\cap (F_{\m}{\times}\R)\,,\qquad (Q^-_\m)':= \Big( Q_{\m}^-+(-16\sqrt{n} k^{-1},16\sqrt{n} k^{-1})^n \Big) \cap (F'_\m{\times}\R)\,,
\end{equation*}
and set
\begin{equation}\label{2004180730}
(u_k)_{\m}:= k\text{-th approximating function for }u_{\m} \text{ on }Q^-_{\m}, \text{ by Theorem~\ref{teo:rough}}
\end{equation}
starting from $u_\m$ defined in \eqref{2004180729} in $(Q^-_\m)'$  as the extension to $\widetilde{\Omega}$ (see \eqref{eq:defappr1} and \eqref{eq:defapprox}). 
Then
\begin{equation}\label{2004180731}
(u_k)_{Q^-}:=\sum_{\m} \chi_{Q^-_{\m}} (u_k)_{\m} \,,
\end{equation}
and, repeating the construction on $Q^+$ to get $(u_k)_{Q^+}$ in $Q^+$,
\begin{equation}\label{2004180744}
(u_k)_{Q}:= \chi_{Q^-}\, (u_k)_{Q^-} + \chi_{Q^+}\, (u_k)_{Q^+}\,.
\end{equation}
We observe that, by Theorem~\ref{teo:rough}, $J_{(u_k)_\m}$ is closed and included in a finite union of boundaries of $n$-dimensional cubes (the bad cubes and the boundary good cubes), and $(u_k)_\m$ is smooth outside its jump set up to the boundary of $Q^-_\m$.
Therefore $J_{(u_k)_Q}$ is closed and included in $\bigcup_\m (J_{(u_k)_\m}\cup \partial Q^-_\m)$
(we will see below that it is enough to take $\Gamma$ and the small sets $J'_\m$, see \eqref{1005181147}, instead of the union of all $\partial Q^-_\m$).

Moreover, $(u_k)_Q\in SBV(Q;\Rn) \cap C^\infty(\ol Q\sm J_{(u_k)_Q};\Rn)\cap W^{m,\infty}(Q\sm J_{(u_k)_Q};\Rn)$ for every $m\in \N$, since this holds separately for each $(u_k)_\m$ up to the boundary of $Q^-_\m$.

Notice that the presence of $\sqrt{n}$ in the sets above is due to the fact that the cubes $Q_j$, oriented with $\nu(x_j)$, are not oriented as the cubes $\qz$ of Theorem~\ref{teo:rough}, which have faces parallel to the axes (see figures at page \pageref{fig}).
Moreover, differently from $Q^-$, the $n$-dimensional measure of $Q^-_{\m}$ vanishes as $k\to\infty$, and at the level of $Q_{\m}^-$ we have to employ the construction of Theorem~\ref{teo:rough} exactly at the scale $k$.
\newline
\\
{\bf Properties of the approximating functions in the cubes.}
For any $\m$, Lemma~\ref{le:Nitsche} gives (as usual we omit the target sets $\Rn$ and $\Mnn$ in the notation for the $L^1$ norm of $u$ and $e(u)$)
\begin{subequations}\label{eqs:2104180829}
\begin{align}
\|u_{\m}\|_{L^1(R'_{\m})}&\leq  C \|u\|_{L^1(R_\m)}\,,\label{2104181100}\\
 \|e(u_{\m})\|_{L^p(R'_{\m})}&\leq  C \|e(u)\|_{L^p(R_\m)}\,,\label{2104181101}\\
 \hn(J_{u_{\m}}\cap R'_{\m}) & \leq C \hn(J_u \cap R_{\m})\,,\label{2104181102}\\
\int \limits_{R'_{\m}} \big|[u_{\m}]\big| \dh &\leq C \int \limits_{R_{\m}} \big|[u]\big| \dh \,.\label{2104181103}
\end{align}
\end{subequations}
By \eqref{2004181946}, \eqref{2004181910}, \eqref{2305181808}, \eqref{2004181940}, and \eqref{eqs:2104180829} 
\begin{equation}\label{2305181921}
\begin{split}
\|u_{\m}-(u_k)_{\m}\|_{L^1(Q^-_{\m})} \leq  \,C k^{-1} |\E u|\big((Q^-_\m)'\sm R'_\m\big)+ C k^{-1/2} \|u\|_{L^1((Q^-_\m)'\sm R'_\m) } + C \|u\|_{L^1(\Omega^2_\m)}\,,
\end{split}
\end{equation}
with \[\mathcal{L}^n(\Omega^2_\m) \leq C \,k^{-1/2} \,\hn\big(J_u \cap \big((Q^-_\m)'\sm R'_\m\big)\big)\,.\]
Summing on $\m$ we get for 
\begin{equation}\label{2104182337}
\widehat{u}_{Q^-}:= \sum_\m \chi_{Q^-_\m} u_\m
\end{equation} 
(notice that the cubes $Q_\m$ overlap at most two times since $\eta_\varepsilon^{-1}$ is larger than $16\sqrt{n}$)
\begin{equation}\label{2204180755}
\begin{split}
\|(\widehat{u}_{Q^-})-(u_k)_{Q^-}\|_{L^1(Q^-)} \leq & \,C k^{-1} |\E u|\big((Q^-+(-t,t)^n)\sm \Gamma\big) + C k^{-1/2} \|u\|_{L^1(Q^-+(-t,t)^n)} \\&+ C \|u\|_{L^1(\Omega^2_{Q^-})}\,,
\end{split}
\end{equation}
with \[\mathcal{L}^n(\Omega^2_{Q^-}) \leq C \,k^{-1/2} \,\hn\big(J_u \cap (Q^-+(-t,t)^n)\sm \Gamma\big)\,.\]
On the other hand, from \eqref{2004180729} and \eqref{2104181100}, 
\begin{equation*}
\|(\widehat{u}_{Q^-}) - u\|_{L^1(Q^-)}\leq C \|u\|_{L^1((Q^-+(-t,t)^n) \cap \{\mathrm{d}(\cdot, \Gamma) < 25\sqrt{n} k^{-1}\})}\,,
\end{equation*}
and then
\begin{equation}\label{2104180838}
\|u-(u_k)_{Q^-}\|_{L^1(Q^-)}\rightarrow 0 \quad\text{as }k\to \infty\,.
\end{equation}

As for $e(u)$, starting from \eqref{2004182216} applied in any $Q^-_\m$, we get
\begin{equation}\label{1105182102}
\begin{split}
\|e((u_k)_{Q^-})\|_{L^p(Q^-)} &\leq \|e(\widehat{u}_{Q^-})\|_{L^p(Q^-+(-t,t)^n)} + C k^{-q} \|e(\widehat{u}_{Q^-})\|_{L^p(Q^-+(-t,t)^n)} \\&\ \ + C\,\theta^q\, \|e(\widehat{u}_{Q^-})\|_{L^p(\Omega^2_{Q^-})}\,,
\end{split}
\end{equation}
and \eqref{2104181101} implies 
\begin{equation}\label{2104180857}
\|e(\widehat{u}_{Q^-}) - e(u)\|_{L^p(Q^-)}\leq C\|e(u)\|_{L^p((Q^-+(-t,t)^n) \cap \{\mathrm{d}(\cdot, \Gamma) < 25\sqrt{n} k^{-1}\})}\,.
\end{equation}
\newline
\\
Let us now estimate the measure and the energy of the jump set of $(u_k)_Q$ in the interior of $Q$.

For any $\m$, \eqref{3rough}, \eqref{4rough} for $Q^-_\m$ and \eqref{2104181102}, \eqref{2104181103} give
\begin{subequations}\label{eqs:2704180915}
\begin{align}
\hn(J_{(u_k)_\m} \cap Q^-_\m)\leq C\,\theta^{-1} &\hn(J_{u_\m} \cap (Q^-_\m)') \leq C\,\theta^{-1} \hn(J_u \cap (Q^-_\m)'\sm R'_\m)\,,\label{2104182322}\\
\int \limits_{J_{(u_k)_\m} \cap Q^-_\m} \hspace{-2em}\big|[(u_k)_\m]\big| \dh \leq C \hspace{-1.5em}&\int \limits_{J_{u_\m} \cap (Q^-_\m)'} \hspace{-1.5em}\big|[u_\m]\big| \dh \leq C \hspace{-2em}\int \limits_{J_u \cap (Q^-_\m)'\sm R'_\m} \hspace{-2em}\big|[u]\big| \dh\,.\label{2104182323}
\end{align}
\end{subequations}
Notice that for the (almost) half cubes $(Q_h^0)^\pm$ we have to consider also the possible jump due to the fact that we have extended $u$ outside $\Omega$ with 0, so we could have created jump on $\partial\Omega\sm \Gamma_h^0$. So the two estimates above include also in the right hand sides the two terms 
\[
C\theta^{-1} \hn(((Q_h^0)^-_\m)' \cap \partial\Omega\sm \Gamma_h^0)\,,\quad \text{and} \qquad C\hspace{-2.5em}\int\limits_{((Q_h^0)^-_\m)' \cap \partial\Omega\sm \Gamma_h^0} \hspace{-2.5em}|\mathrm{tr}_\Omega u| \dh\,,
\]
respectively.

Let us examine the jump for $(u_k)_{Q^-}$ created on the common boundaries between two sets $Q^-_\m$, $Q^-_{\m'}$, namely between two sets $Q^-_\m$ and $Q^-_{\m\pm e_i}$ for $i=1,\dots,{n-1}$, both inside $Q^-$. To fix the ideas let us take $\m$ and consider $Q^-_\m$ and $Q^-_{\m+ e_i}$. Notice that 
\[
|m_n(\m)-m_n(\m+e_i)| \leq 1/2\,,
\]
where $m_n(\m)k^{-1}$ and $m_n(\m+e_i)k^{-1}$ are the ``heights'' corresponding to $Q^-_\m$ and $Q^-_{\m+ e_i}$, see \eqref{1904182312}. This means that, for $m_n=m_n(\m)$,
\begin{equation*}
\begin{split}
u_\m&=u \text{ in } F'_{\m}{\times}(-\varrho-16\sqrt{n} k^{-1}, (m_n-1/2) k^{-1})\,,\\ u_{\m+e_i}&=u \text{ in } F'_{\m+e_i}{\times}(-\varrho-16\sqrt{n} k^{-1}, (m_n-1/2) k^{-1})\,.
\end{split}
\end{equation*}
By construction of $(u_k)_\m$ (see \eqref{eq:defappr1} and \eqref{eq:defapprox}) we have that 
\begin{equation}\label{2104181251}
(\widetilde{u}_k)_\m=(\widetilde{u}_k)_{\m+e_i}\, \text{ in }(F'_\m \cap F'_{\m+e_i}) {\times} (-\varrho-16\sqrt{n} k^{-1}, (m_n-(4\sqrt{n}+1/2) k^{-1})\,,
\end{equation}
and 
\[(u_k)_\m=(u_k)_{\m+e_i}\, \text{ in }(F'_\m \cap F'_{\m+e_i}) {\times} (-\varrho-16\sqrt{n} k^{-1}, (m_n-(8\sqrt{n}+1/2) k^{-1})\,,\]
since, if $x \in \qz$, $(\widetilde{u}_k)_\m(x)$ and $(u_k)_\m(x)$ depend only on $u_\m$ in $\tqz$ and $\Qz$, respectively (see figure on the right at page \pageref{fig}).

Setting
\begin{equation}\label{1005181147}
J_{(u_k)_{Q^-}}\cap \partial Q^-_\m \cap \partial Q^-_{\m+e_i}=:J'_\m\,,
\end{equation}
it follows that 
\begin{equation}\label{2104181328}
J'_\m \subset (\partial F_\m \cap \partial F_{\m+e_i}){\times}\big( (m_n-(8\sqrt{n}+1/2), m_n+1)k^{-1} \big)
\end{equation}
and thus
\begin{equation}\label{2104182149}
\hn(J'_\m)\leq C \eta_\varepsilon^{-(n-2)} k^{-(n-1)} \,.
\end{equation}
Summing up over all the faces of $Q^-_\m$ in the directions $e_1,\dots, e_{n-1}$ we get
\begin{equation}\label{2604182000}
\hn(J_{(u_k)_{Q^-}}\cap \partial Q^-_\m)\leq C \eta_\varepsilon^{-(n-2)} k^{-(n-1)} \,,
\end{equation}
and summing up over $\m$ gives (see \eqref{2904181212})  
\begin{equation}\label{2104181234}
\hn\Big(J_{(u_k)_{Q^-}}\cap \bigcup_\m \partial Q^-_\m \sm \partial Q^-\Big)\leq C \eta_\varepsilon \varrho^{n-1} \,.
\end{equation}

In order to estimate the amplitude of the jump, we start from the jump on $J'_\m$. For every $x\in J'_\m$ we may have four cases, depending if $x \in (\widetilde{\Omega}_b^k)_\m$, $x \in (\widetilde{\Omega}_b^k)_{\m+e_i}$, or not, where $(\widetilde{\Omega}_b^k)_\m$ is the set of (neighbourhoods of) bad cubes corresponding to $(Q^-_\m)'$, see \eqref{1204181406}.
By construction of $(u_k)_\m$ it follows that
\begin{equation*}
[(u_k)_{Q^-}]= \varphi_k \ast \big( (\widetilde{u}_k)_\m -  (\widetilde{u}_k)_{\m+e_i} \big)\quad\text{in }J'_\m\sm\big((\widetilde{\Omega}_b^k)_\m \cup(\widetilde{\Omega}_b^k)_{\m+e_i}\big)\,,
\end{equation*}
so, for every $x$ in the set above,
\begin{equation}\label{2104182147}
\big|[(u_k)_{Q^-}] \big| \leq \|\varphi\|_{L^\infty(B(0,1))} k^n \|(\widetilde{u}_k)_\m -  (\widetilde{u}_k)_{\m+e_i} \|_{L^1(J'_\m+ B(0,k^{-1}))}\,.
\end{equation}
We claim that (see \eqref{2104181322} for the definition of $R_\m$) 
\begin{equation}\label{2104181313}
\|(\widetilde{u}_k)_\m -  (\widetilde{u}_k)_{\m+e_i} \|_{L^1(J'_\m + B(0,k^{-1}))} \leq  C k^{-1} |\E u|(R_\m \cap R_{\m+e_i})\,.
\end{equation}
We have (in the following we consider any enumeration $(z_j)_j$ for the nodes $z$, and denote the cubes corresponding to $z_j$ by $q_j$, $\tilde{q}_j$, $Q_j$, $\widetilde{Q}_j$, in no context with the  cubes of scale $\varrho_j$ introduced in \eqref{eqs:2405171202})
\begin{equation}\label{2604181109}
\begin{split}
\|(&\widetilde{u}_k)_\m -  (\widetilde{u}_k)_{\m+e_i} \|_{L^1(J'_\m + B(0,k^{-1}))} \leq \hspace{-2em}\sum_{q_j\cap (J'_\m + B(0,k^{-1}))\neq \emptyset} \hspace{-2em}\|(\widetilde{u}_k)_\m -  (\widetilde{u}_k)_{\m+e_i} \|_{L^1(q_j)}\\& 
\leq \hspace{-2em}\sum_{q_j\cap (J'_\m + B(0,k^{-1}))\neq \emptyset} \hspace{-2em}\|(\tilde{a}_j)_{\m} - (\tilde{a}_{j})_{\m+e_i}\|_{L^1(q_j)} + C k^{-1} \hspace{-2em}\sum_{q_j\cap (J'_\m + B(0,k^{-1}))\neq \emptyset} \hspace{-2em}(|\E(u_\m)|(\widetilde{Q}_j) + |\E (u_{\m+e_i})|(\widetilde{Q}_j))
\end{split}
\end{equation}
where $(\tilde{a}_j)_\m$ affine with $e((\tilde{a}_j)_\m)=0$ and
\begin{equation*}
\|u_\m - (\tilde{a}_j)_\m\|_{L^1(\tilde{q}_j)} \leq C k^{-1}|\E u_\m|(\tilde{q}_j)\,.
\end{equation*}
The second inequality in \eqref{2604181109} comes from (let $(\omega^k)_\m$ be the exceptional set corresponding to $u_\m$, see Theorem~\ref{teo:rough}) 
\begin{equation*}
\|(\widetilde{u}_k)_\m - (\tilde{a}_j)_\m\|_{L^1(q_j\sm (\omega^k)_m)}=\|u_\m- (\tilde{a}_j)_\m\|_{L^1(q_j\sm (\omega^k)_m)} \leq C k^{-1}|\E u_\m|(\tilde{q}_j)\,,
\end{equation*}
and the fact that, recalling \eqref{1204182345},
\begin{equation*}
\|(\widetilde{u}_k)_\m - (\tilde{a}_j)_\m\|_{L^1(q_j\cap (\omega^k)_m)}\leq C k^{-1} \big(|\E u_\m|(q_j)+ \|e(u_\m)\|_{L^1(\widetilde{Q}_j)}\big)\,,
\end{equation*}
the same being true for $\m + e_i$ in place of $\m$.

We now estimate $\|(\tilde{a}_j)_{\m} - (\tilde{a}_{j})_{\m+e_i}\|_{L^1(q_j)}$ for $q_j\cap (J'_\m + B(0,k^{-1}))\neq \emptyset$ in \eqref{2604181109}. We remark that
\[\mathcal{L}^n(\tilde{q}_j\cap R_\m \cap R_{\m+e_i})/\mathcal{L}^n(\tilde{q}_j)\geq C_0 >0\,,\]
with $C_0$ depending only on $n$.
Thus
\begin{equation}\label{1805180205}
\begin{split}
\|(\tilde{a}_j)_{\m} - (\tilde{a}_{j})_{\m+e_i}\|_{L^1(q_j)}&\leq C\|(\tilde{a}_j)_{\m} - (\tilde{a}_{j})_{\m+e_i}\|_{L^1(\tilde{q_j}\cap R_\m\cap R_{\m+e_i})}\\
&\leq C k^{-1}(|\E(u_\m)|(\widetilde{Q}_j) + |\E (u_{\m+e_i})|(\widetilde{Q}_j))\,,
\end{split}
\end{equation}
since $(\tilde{a}_j)_{\m} - (\tilde{a}_{j})_{\m+e_i}$ is an affine function and $u_\m=u_{\m+e_i}=u$ in $R_\m\cap R_{\m+e_i}$. (The constant in the first inequality above depends on $C_0$.)
%
Therefore \eqref{2104181313} is proven, recalling also \eqref{eqs:2104180829}.

Consider now the case when $x \in J'_\m \cap (\widetilde{\Omega}_b^k)_\m $. To fix the ideas assume that $x\in q_j$ (in the open cube). So (recall \eqref{eq:defapprox})
\begin{equation*}
(u_k)_\m(x)=(\tilde{a}_j)_\m(x)\,,\quad\text{with } \|u_\m - (\tilde{a}_j)_\m\|_{L^1(\tilde{q}_j)} \leq C k^{-1}|\E u_\m|(\tilde{q}_j)\,.
\end{equation*} 
If $x\notin (\widetilde{\Omega}_b^k)_{\m+e_i}$, $(u_k)_{\m+e_i}(x)=\varphi_k\ast (\widetilde{u}_k)_{\m+e_i}$, so 
\[
[(u_k)_{Q^-}](x)= \varphi_k \ast \big((\widetilde{u}_k)_{\m+e_i} - (\tilde{a}_j)_\m \big)(x)\,.\]
Now
\begin{equation*}
\begin{split}
\|(&\widetilde{u}_k)_{\m+e_i} - (\tilde{a}_j)_\m \|_{L^1(B(x,k^{-1}))} \leq \|(\widetilde{u}_k)_{\m+e_i} - (\tilde{a}_j)_\m \|_{L^1(\tilde{q}_j)}\\&\leq \|(\tilde{a}_j)_{\m+e_i} - (\tilde{a}_j)_\m \|_{L^1(\tilde{q}_j)} 
+ C k^{-1}|\E u_{\m+e_i}|(\tilde{q}_j)
\leq C k^{-1} (|\E(u_\m)|(\widetilde{Q}_j) + |\E (u_{\m+e_i})|(\widetilde{Q}_j))\,,
\end{split}
\end{equation*}
arguing as done for \eqref{2604181109} and \eqref{1805180205}.
In the 
same way one deals with the case $x \in J'_\m \cap (\widetilde{\Omega}_b^k)_{\m+e_i}\sm (\widetilde{\Omega}_b^k)_\m$. The last case is $x\in  J'_\m\cap (\widetilde{\Omega}_b^k)_{\m+e_i}\cap (\widetilde{\Omega}_b^k)_\m$: now
\begin{equation*}
\begin{split}
[(u_k)_{Q^-}](x)=|(\tilde{a}_j)_{\m+e_i} - (\tilde{a}_j)_\m|(x)\,.
\end{split}
\end{equation*}

We now put together the different cases, deducing that 
\begin{equation*}
\big|[(u_k)_{Q^-}]\big| \leq C k^{n-1} |\E u|(R_\m \cap R_{\m+e_i}) \quad\text{in }J'_\m\,,
\end{equation*} 
so that \eqref{2104182149} gives, integrating over $J'_\m$, that
\begin{equation}\label{2604182001}
\int \limits_{J'_\m} \big|[(u_k)_{Q^-}]\big| \dh \leq C \eta_\varepsilon^{-(n-2)}|\E u|(R_\m \cap R_{\m+e_i})\,.
\end{equation} 
Since in the estimates are employed the cubes $\widetilde{Q}_j$, with sidelength $16k^{-1}$, we look possibly at height $16\sqrt{n}k^{-1}$ below $J'_\m$, which is distant less than $9\sqrt{n}k^{-1}$ from $\Gamma$. This motivates the choice of the constant 25 in the definition of $R_\m$.

Summing up over all the faces of $Q^-_\m$ in the directions $e_1,\dots, e_{n-1}$ and over $\m$ (observe that $R_\m \cap R_{\m\pm e_i}$ overlap each other at most 2 times, over $i$ and $\m$) we deduce
\begin{equation}\label{2204180924}
\int \limits_{J_{(u_k)_{Q^-}}\cap \bigcup_\m \partial Q^-_\m \sm \partial Q^-} \big|[(u_k)_{Q^-}]\big| \dh \leq C \eta_\varepsilon^{-(n-2)}|\E u|(\{\mathrm{d}(\cdot, \Gamma) < 25\sqrt{n} k^{-1}\} \sm \ol{Q^+})\,.
\end{equation}

Let us now consider the jump of $(u_k)_Q$ on $\Gamma$, by looking separately at the traces of $u-(u_k)_{Q^\pm}$ on the two sides of $\Gamma$. 
We have ($\tr^-$ denotes the trace on $\Gamma$ from $Q^-$)
\begin{equation*}
\int \limits_{\Gamma\cap Q^-} \tr^-(u-(u_k)_{Q^-}) \dh = \int \limits_{\Gamma\cap Q^-} \tr^-(u-\widehat{u}_{Q^-}) \dh + \int \limits_{\Gamma\cap Q^-} \tr^-((\widehat{u}_{Q^-})-(u_k)_{Q^-}) \dh
\end{equation*}
where $\widehat{u}_{Q^-}$ has been introduced in \eqref{2104182337}.
By definition \eqref{2004180729} of $u_\m$ one has
\begin{equation*}
\int \limits_{\Gamma \cap Q^-_\m} |\tr^-(u-u_\m)|\dh \leq C|\E(u-u_\m)|(F'_\m{\times}\{\mathrm{d}(\cdot, \Gamma)< 2 k^{-1}\})\leq C |\E u|(R_\m)\,, 
\end{equation*}
and, summing up over $\m$,
\begin{equation}\label{2204180757}
\int \limits_{\Gamma\cap Q^-} \tr^-(u-\widehat{u}_{Q^-}) \dh \leq C |\E u|\big((Q^-+(-t,t)^n) \cap \{\mathrm{d}(\cdot, \Gamma)< 25 \sqrt{n} k^{-1}\}\sm \ol{Q^+}\big)\,.
\end{equation}
Moreover, arguing as in \cite[Theorem~3.2, Steps 1 and 4]{Bab15} (see also the proof of \cite[Theorem~1.1, property (1.1d)]{CC17}), we get
\begin{equation}\label{2204180021}
\begin{split}
\int \limits_{\Gamma\cap Q^-} |\tr^-((\widehat{u}_{Q^-})- & (u_k)_{Q^-})| \dh\leq \frac{C}{t}\|(\widehat{u}_{Q^-})-(u_k)_{Q^-}\|_{L^1((Q^-+(-t,t)^n)\sm \ol{Q^+})} \\&+ |\E\big((\widehat{u}_{Q^-})-(u_k)_{Q^-}\big)| ((Q^-+(-t,t)^n)\cap\{\mathrm{d}(\cdot, \Gamma)<t \} \sm \ol{Q^+})\,.
\end{split}
\end{equation}
Collecting \eqref{2204180757} and \eqref{2204180021} we estimate $\tr^-(u-(u_k)_{Q^-})$ on $\Gamma\cap Q^-$. Arguing in the same way for the positive trace (namely, that corresponding to $Q^+$)
and adding the two, we obtain
\begin{equation}\label{2204180800}
\begin{split}
\int \limits_{\Gamma\cap Q} \big|[u]-[(u_k)_{Q}]\big| \dh &\leq C |\E u|\big( (Q+(-t,t)^n)  \cap \{\mathrm{d}(\cdot, \Gamma)<t \} \sm \Gamma\big) \\&+ \frac{C}{t} \|\widehat{u}_{Q}-(u_k)_{Q}\|_{L^1(Q+(-t,t)^n)}\,,
\end{split}
\end{equation}
setting $\widehat{u}_{Q}:=\chi_{Q^-}\widehat{u}_{Q^-} + \chi_{Q^+}\widehat{u}_{Q^+}$ (and $\widehat{u}_{Q^+}$ defined in analogy to $\widehat{u}_{Q^-}$). 
If we are in a boundary cube $Q_h^0$, we consider $u$ extended with 0 outside $\Omega$, so that on $\partial\Omega$ we replace $[u]$ with $\tr_\Omega u$ also in the right hand side of \eqref{2204180800}, in the evaluation of $|\E u|$. 
\newline
\\
{\bf The approximating functions.}
We consider $\widetilde{B}_0:=B_0+(-t,t)^n$ and we denote $(u_k)_{B_0}$ the $k$-th approximating function for $u$ given by Theorem~\ref{teo:rough} in correspondence to $B_0$, starting from the 
the extension of $u$, with value 0 outside $\Omega$, in $\widetilde{B}_0$. Notice that 
\eqref{3rough} and \eqref{4rough} give
\begin{subequations}
\begin{align}
\hn(J_{(u_k)_{B_0}})& \leq C \,\theta^{-1} \hn\big((J_u \cap \widetilde{B}_0) \cup (\partial\Omega\cap \widetilde{B}_0)\big)\,,\label{2704180939}\\
\limsup_{k\to \infty} \int\limits_{J_{(u_k)_{B_0}}} \hspace{-1em}\big|[(u_k)_{B_0}]\big| \dh &\leq C \Big(\int\limits_{J_u \cap \widetilde{B}_0} \hspace{-1em}\big|[u]\big| \dh + \int \limits_{\partial\Omega\cap \widetilde{B}_0} \hspace{-1em}|\mathrm{tr}_\Omega u|\dh \Big)\,, \label{2704180940}
\end{align}
\end{subequations}
respectively.
Then we define the global $k$-th approximating function
\begin{equation}\label{eq:defApprFct}
u_k:=\chi_{B_0} (u_k)_{B_0} + \sum_{j=1}^{\ol \jmath} \chi_{Q_j} (u_k)_{Q_j} + \sum_{h=1}^{\ol h} \chi_{Q_h^0 \cap \Omega} (u_k)_{Q_h^0}\,,
\end{equation}
where $(u_k)_{Q}$
are introduced in \eqref{2004180744}.

Notice that the functions $(u_k)_Q$ and $(u_k)_{B^0}$ are smooth up to the boundaries of their domains, outside their jump sets, which are closed and included in a finite union of $C^1$ hypersurfaces (see the discussion below \eqref{2004180744}).

Then $u_k\in SBV(\Omega;\Rn)\cap C^\infty(\ol\Omega\sm J_{u_k};\Rn) \cap W^{m,\infty}(\Omega\sm J_{u_k};\Rn)$ for every $m\in \N$, $J_{u_k}$ is closed and 
\begin{equation*}
J_{u_k}\subset J_{(u_k)_{B^0}} \cup \bigcup_Q (J_{(u_k)_Q} \cup \partial Q)
\end{equation*}
(where $Q$ stands for all the $Q_j$ and $Q_h^0$) which is a finite union of $C^1$ hypersurfaces (we will see below that it is enough to take just a little part of $\partial Q$, see \eqref{2204180915}).

By definitions \eqref{2004180744} and \eqref{2104180838} we have that $(u_k)_{Q} \to u$ in $L^1(Q;\Rn)$ for every $Q$.
Moreover, $(u_k)_{B_0} \to u$ in $L^1(B_0;\Rn)$ by \eqref{1rough}, so \eqref{eq:defApprFct} and \eqref{eq:defb0} imply
\begin{equation}\label{2004181411}
u_k \to u \quad\text{in }L^1(\Omega;\Rn)\,.
\end{equation}
We can argue very similarly to prove \eqref{5main}, starting from \eqref{5rough} applied in each $Q^-_\m$ (this gives the analogous of \eqref{2104181100}, then we follow the argument for \eqref{2104180838}).

Putting together \eqref{1105182102} and \eqref{2104180857} raised to the power $p$ for any $Q$, and collecting with \eqref{2rough} for $B_0$ (again raised to the $p$) we obtain immediately
\begin{equation}\label{2104180906}
\limsup_{k\to \infty} \|e(u_k)\|_{L^p(\Omega;\Mnn)}\leq \|e(u)\|_{L^p(\Omega;\Mnn)}\,.
\end{equation}

In order to treat the jump set, notice that we have still to estimate the jump of $u_k$ on $\bigcup_j \partial Q_j \cup \bigcup_h (\partial Q_h^0 \cap\Omega)$. To do so, we may closely follow what done for the jump on $\partial Q^-_\m$: the only difference is that now  we have in $B_0$ the rough approximation of $u$, without any extension in the spirit of Lemma~\ref{le:Nitsche}. Then, if we have two parallelepipeds $Q_\m\subset Q$ and $Q_{\m+e_i}\cap B_0\neq \emptyset$, we consider in \eqref{1904182312}
\begin{equation}
u_{\m+e_i}=u \,\text{ in } F'_{\m+e_i}{\times}(-\varrho-16\sqrt{n}k^{-1}, (m_n + 25 \sqrt{n})k^{-1})\,.
\end{equation}
Differently from before, now $|\E(u_\m)|(\widetilde{Q}_j) + |\E (u_{\m+e_i})|(\widetilde{Q}_j)$,
entering for instance in \eqref{2604181109}, is estimated by $|\E u|(R_\m\cup (R_{\m+e_i} \cup R'_{\m+e_i}))$, see \eqref{2104181322} for the definition of $R'_\m$. For this reason, for the analogue of \eqref{2604182001} we get
\begin{equation}\label{2204180925}
\begin{split}
\int \limits_{J_{u_k}\cap \partial Q_j} \big|[u_k]\big| \dh &\leq C \eta_\varepsilon^{-(n-2)} \Big[|\E u|\big( (Q_j + (-t,t)^n) \cap \widetilde{B}_0 \cap \{\mathrm{d}(\cdot, \Gamma_j) < 25\sqrt{n} k^{-1}\} \sm \Gamma_j \big) \\& \hspace{1em} + |\E u|(\Gamma_j \cap \{\mathrm{d}(\cdot, \partial Q_j)<32 \sqrt{n} k^{-1}\}) \Big]\,,
\end{split}
\end{equation}
and the same for $Q_h^0$, $\Gamma_h^0$ in place of $Q_j$, $\Gamma_j$. Notice that we have an additional term with respect to \eqref{2604182001}, which vanishes as $k$ tends to $\infty$, since $|E^j u|$ is evaluated on a subset of $\Gamma_j$ whose $\hn$ measure vanishes in $k$.

As done for \eqref{2604182000}, one deduces
\begin{equation}\label{2204180915}
\hn(J_{u_k}\cap \partial Q_j) \leq C \eta_\varepsilon^{-(n-2)} k^{-(n-1)}\,,\quad\hn(J_{u_k}\cap \partial Q_h^0 \cap \Omega) \leq C \eta_\varepsilon^{-(n-2)} k^{-(n-1)} \,.
\end{equation}

We have
\begin{equation*}
J_{u_k} \subset (J_{u_k}\cap B_0) \cup \bigcup_{j=1}^{\ol \jmath} (J_{u_k}  \cap \ol Q_j \sm \Gamma_j) \cup \widehat{\Gamma} \cup \bigcup_{h=1}^{\ol h} (J_{u_k} \cap \ol{Q}_h^0\cap \Omega)\,.
\end{equation*}
Moreover, we may assume that $\widehat{\Gamma} \subset J_{u_k} $, since there are arbitrarily small $a>0$ with $\hn( \widehat{\Gamma} \cap \{[u_k]=a\})=0$, and then we can add to $u_k$ a perturbation with arbitrarily small $W^{1,\infty}(\Omega\sm \widehat{\Gamma})$ norm, having jump of class $C^1$ on $\widehat{\Gamma}$ and equal to $a$ on an arbitrarily large subset of $\widehat{\Gamma}$ (see also \cite[Lemmas~4.1, 4.3]{DPFusPra17}).
Therefore we may assume that
\begin{equation}\label{2204181158}
J_{u_k} \triangle J_u \subset  (J_{u_k} \sm \widehat{\Gamma}) \cup (J_u \triangle \widehat{\Gamma}) \,.
\end{equation}
Collecting \eqref{2104182322}, \eqref{2104181234}, \eqref{2204180915}, \eqref{2704180939}, and recalling \eqref{1305171150}, we deduce
\begin{equation*}
\hn\big(J_{u_k}\sm \widehat{\Gamma}) \leq  C\, \theta^{-1} \big(\hn(J_u \sm \widehat{\Gamma}) + \hn(\partial\Omega\sm \widehat{\Gamma}_{\partial\Omega})\big) + C\, \big(\hn(J_u) + \hn(\partial\Omega)\big) \, \eta_\varepsilon\,.
\end{equation*}
By \eqref{2204181223} it then follows that
\begin{equation*}\label{2204181152}
\hn(J_{u_k}\triangle J_u) \leq C \, \theta^{-1} \varepsilon + C\, \eta_\varepsilon\,.
\end{equation*}
As $\varepsilon$ is arbitrary and $\lim_{\varepsilon\to 0}\eta_\varepsilon=0$, we conclude
\begin{equation}\label{2904181256}
\lim_{k\to \infty} \hn(J_{u_k}\triangle J_u)=0\,.
\end{equation}

The combination of \eqref{2104182323}, \eqref{2204180924}, \eqref{2204180925}, and \eqref{2704180940} gives
\begin{equation}\label{2204181346}
\begin{split}
\int \limits_{J_{u_k} \sm \widehat{\Gamma}} \big| [u_k] \big| \dh &\leq C \big(1+\eta_\varepsilon^{-(n-2)}\big) \int \limits_{(J_u \sm  \widehat{\Gamma})\cup \widehat{\Gamma}_k} \big|[u]\big| \dh + C\, \eta_\varepsilon^{-(n-2)} \|e(u)\|_{L^1(\{ \mathrm{d}(\cdot, \widehat{\Gamma}\cup \partial \Omega) < 25\sqrt{n} k^{-1} \}) }
\\ & \leq C\, \eta_\varepsilon + C\, \eta_\varepsilon^{-(n-2)} \Big(\|e(u)\|_{L^1(\{ \mathrm{d}(\cdot, \widehat{\Gamma}\cup \partial \Omega) < 25\sqrt{n} k^{-1} \}) }+ \int\limits_{\widehat{\Gamma}_k} \big|[u]\big| \dh\Big)\,,
\end{split}
\end{equation}
letting $\widehat{\Gamma}_k:=\bigcup_j (\Gamma_j \cap \{\mathrm{d}(\cdot, \partial Q_j)<C k^{-1}\}) \cup \bigcup_h (\Gamma_h^0 \cap \{\mathrm{d}(\cdot, \partial Q_h^0)<C k^{-1}\})$, and recalling the definition \eqref{2204181256} of $\eta_\varepsilon$. Notice that in the first inequality in \eqref{2204181346} we should have written all the term in \eqref{2204181256}, which is nothing but the jump part of the extension of $u$ with 0 outside $\Omega$ (see also the remark below  
\eqref{eqs:2704180915}).

Summing up \eqref{2204180800} for $j=1,\dots, \ol \jmath$ and employing \eqref{2204180755} we get 
\begin{equation}\label{2204181350}
\begin{split}
\int \limits_{\widehat{\Gamma}} \big|[u]-[u_k]\big| \dh & \leq C\, \|e(u)\|_{L^1(\{\mathrm{d}(\cdot, \widehat{\Gamma})<t\})} + C \int \limits_{J_u \sm \widehat{\Gamma}} \big|[u]\big| \dh \\& + \frac{C}{t}\Big( k^{-1} |\E u| (\Omega\sm \widehat{\Gamma}) + k^{-1/2} \|u\|_{L^1(\Omega)} + \|u\|_{L^1(\Omega^2)}  \Big)\,,
\end{split}
\end{equation}
with
\begin{equation*}
\mathcal{L}^n(\Omega^2)\leq C\, k^{-1/2} \,\hn(J_u \sm \widehat{\Gamma})\,.
\end{equation*}
Since
\begin{equation*}
\int \limits_{J_u \cup J_{u_k}} \big|[u]-[u_k]\big| \dh \leq \int \limits_{\widehat{\Gamma}} \big|[u]-[u_k]\big| \dh + \int \limits_{J_{u_k} \sm \widehat{\Gamma}} \big| [u_k] \big| \dh + \int \limits_{J_u \sm  \widehat{\Gamma}} \big|[u]\big| \dh\,,
\end{equation*}
we conclude that
\begin{equation}\label{2704181126}
\lim_{k\to \infty} \int \limits_{J_u \cup J_{u_k}} \big|[u]-[u_k]\big| \dh = 0
\end{equation}
collecting \eqref{2204181346} and \eqref{2204181350} and sending $k\to \infty$, $t\to 0$, and $\varepsilon \to 0$ in this order.

At this stage we can say that $u_k$ is a sequence bounded in $BD(\Omega)$, converging to $u$ in $L^1(\Omega;\Rn)$ (see \eqref{2004181411}). Therefore, by \cite[Theorem~1.1]{BelCosDM98} and recalling \eqref{2104180906}, this gives 
\begin{equation}\label{2904181259}
\lim_{k\to\infty}\|e(u_k)-e(u)\|_{L^p(\Omega;\Mnn)}=0\,.
\end{equation}
By \eqref{2704181126} we have that 
\begin{equation*}
\lim_{k\to \infty}|\E^j (u-u_k)|(\Omega)=0\,,
\end{equation*}
and then $\|u_k - u\|_{BD(\Omega)}\to 0$. Recalling \eqref{2904181256} and \eqref{2904181259} we conclude the proof.
\end{proof}
\begin{remark}
Looking at the proof of Theorem~\ref{teo:density}, one needs that $\Omega$ has finite perimeter, that there is a suitable notion of trace on $\partial\Omega$, and that the function $u$ considered has trace integrable on $\partial\Omega$.
This would permit to weaken the assumption that $\Omega$ is a bounded Lipschitz domain.
\end{remark}
%
%

\section{Proof of the other density theorems}\label{sec:other}

In this section we discuss two further density results for functions in $SBD(\Omega)$ and in $SBD^p_\infty(\Omega)$ in the spirit of \cite{DPFusPra17}. The space $SBD^p_\infty(\Omega)$ consists of all functions $u\in SBD(\Omega)$ with $e(u)\in L^p(\Omega;\Mnn)$, and without any constraint on $\hn(J_u)$ (see Section~\ref{Sec1}). These results are obtained by corresponding modifications of the rough approximation result Theorem~\ref{teo:rough}, that permit then to follow the strategy of Theorem~\ref{teo:density}.

We assume that $\Omega \subset \Rn$ is a Lipschitz domain. As above, this may be avoided by requiring that $\Omega$ has finite perimeter, that there is a suitable notion of trace on $\partial\Omega$, and that the function $u$ considered has trace integrable on $\partial\Omega$.

The first part of the proof is common for the two results. Since now  $\hn(J_u)$ may be infinite, but we are interested in the approximation in energy, we consider for a fixed $\varepsilon>0$ a set $\widetilde{\Gamma}_\varepsilon \subset J_u$, with $\hn(\widetilde{\Gamma}_\varepsilon)<\infty$, such that
\begin{equation}\label{1105182147}
\int \limits_{J_u\sm \widetilde{\Gamma}_\varepsilon} \big| [u] \big| \dh <\varepsilon\,.
\end{equation} 
This follows from the fact that $[u]\in L^1(J_u;\Rn)$.
Then we employ the approximation procedure at the beginning of proof of Theorem~\ref{teo:density} to $\widetilde{\Gamma}_\varepsilon$ in place of $J_u$ (and to $\partial\Omega$ as before), obtaining a finite family of pairwise disjoint closed cubes $(\ol{Q_j})_{j=1}^{\ol{\jmath}}\subset \Omega$ satisfying the same properties as before (we keep the same notation), with $J_u$ replaced by $\widetilde{\Gamma}_\varepsilon$ (also in \eqref{2204181223}). In particular
\begin{equation}\label{1205182222}
\lim_{\varepsilon\to 0} \int \limits_{J_u\sm \widehat{\Gamma}} \big| [u] \big| \dh = 0\,.
\end{equation}
The definition of $\eta_\varepsilon$ in \eqref{2204181256} remains the same, and $\eta_\varepsilon$ is still vanishing as $\varepsilon\to 0$ thanks to \eqref{1105182147}. 
Notice that we keep the same notation of Theorem~\ref{teo:density}, for instance for the (almost) parallelepipeds $Q^-_\m$ and for the convolution kernel $\varphi_k$.
\begin{proof}[Proof of Theorem~\ref{teo:TeorA}]
Since we are now proving an estimate which is linear both in $e(u)$ and in $\E^j u$, the construction for Theorem~\ref{teo:rough} may be replaced simply by the convolution with $\varphi_k$. 
Indeed for every $v\in SBD(\widetilde{U})$ with $\ol U \subset \widetilde{U}$ we have that, for $k$ large enough, $v_k:=v\ast \varphi_k$ is in $C^\infty(\ol U;\Rn)$ and satisfies
\begin{equation}\label{1105182048}
\int \limits_U |e(v_k)|\dx \leq |\E v|(U+B(0,k^{-1}))\,.
\end{equation}
So we keep all as in Theorem~\ref{teo:density} except for the definition of $(u_k)_\m$ in $Q^-_\m$, given in \eqref{2004180730}: now
\begin{equation}
(u_k)_\m:= u_\m \ast \varphi_k\,,
\end{equation}  
where $u_\m$ is still defined as in \eqref{2004180729} and \eqref{2104181322} (notice that now we could have taken also $R_\m$ of height $\sqrt{n}k^{-1}$ instead of $25 \sqrt{n}k^{-1}$, but we prefer to keep the same notation).

Similarly to before, 
we have that 
\begin{equation}\label{1105182127}
\|u_{\m}-(u_k)_{\m}\|_{L^1(Q^-_{\m})} \leq  \,C k^{-1} |\E u|\big((Q^-_\m)'\sm R'_\m\big)\,,
\end{equation}
while, in place of \eqref{1105182102},
\begin{equation}\label{1105182128}
\|e((u_k)_{Q^-})\|_{L^1(Q^-)} \leq \|e(\widehat{u}_{Q^-})\|_{L^1(Q^-+(-t,t)^n)} + 2\,|\E u|(Q^-\sm \Gamma)\,.
\end{equation}
Notice that \eqref{2104181101} and \eqref{2104180857} hold with the norm $L^1$ instead of the norm $L^p$.

Since now we have not distinguished the cubes in bad and good ones, we have no jump in (the open set) $Q^-_\m$, so \eqref{eqs:2704180915} are useless, and in order to estimate $[u_k]$ on $J'_\m$ (see \eqref{1005181147}) we have only one case, corresponding to the estimate \eqref{2104181313}, which is still true. Also \eqref{2204180800} holds as before.

The approximating functions $u_k$ are defined as in \eqref{eq:defApprFct}, with $(u_k)_{B_0}$ still obtained by convolution between $\varphi_k$ and the function $u$ in $\widetilde{B}_0$, extended with 0 outside $\Omega$. 

Now \eqref{2004181411} and \eqref{2104180906} (with the norm $L^1$ instead of $L^p$) follow from \eqref{1105182127} and \eqref{1105182128}, respectively, 
employing also \eqref{1205182222}.

By (the anologues of) \eqref{2204180924} and \eqref{2204180925} we deduce \eqref{2204181346}, recalling also the definition of $\eta_\varepsilon$ \eqref{2204181256}.

Putting together \eqref{2104181234} and \eqref{2204180915} (that hold also in the present setting) we obtain
\begin{equation}\label{1105182206}
\hn(J_{u_k}\sm \widehat{\Gamma})\leq C \big(\hn(J_u)+\hn(\partial\Omega) \big) \, \eta_\varepsilon\,.
\end{equation}

Moreover, \eqref{2704181126} follows as before from \eqref{2204181346}, that still holds, and \eqref{2204181350}, which is slightly modified since now combines  \eqref{2204180800} and \eqref{1105182127} (instead of \eqref{2305181921}). Since $u_k$ is bounded in $BD(\Omega)$, then 
\eqref{2004181411}, \eqref{2104180906}, \eqref{2704181126}, and \eqref{1105182206} give \eqref{eqs:mainTeorA}.

It lasts only to prove that $J_{u_k}$ is, up to a negligible set,
 a finite union of \emph{pairwise disjoint} compact $C^1$ hypersurfaces contained in $\Omega$. To do so, notice that
\begin{equation}\label{1105182200}
J_{u_k}\subset \widehat{\Gamma} \cup \bigcup_Q \bigcup_\m J'_\m \subset \subset \Omega\,,
\end{equation}
because there is not the jump due to bad cubes and boundary good cubes in any $Q^-_\m$ and in $B_0$. 
Since $J'_\m$ are in a finite number and transversal to $\widehat{\Gamma}$, we have that $\widehat{\Gamma} \cap \bigcup_Q \bigcup_\m J'_\m$ consists in a finite number of $n{-}2$ dimensional manifolds, with finite $\mathcal{H}^{n{-}2}$ measure. Therefore we may follow the capacitary argument by Cortesani in \cite[Corollary~3.11]{Cor97}, replacing the jump in a small neighbourhood of $\widehat{\Gamma} \cap \bigcup_Q \bigcup_\m J'_\m$ by an $H^1$ transition with arbitrary small $H^1$ norm (this is possible since the capacitary argument is applied to $u_k\in L^\infty(\Omega;\Rn)$ and since the $2$-capacity of $\widehat{\Gamma} \cap \bigcup_Q \bigcup_\m J'_\m$ is 0, because it has finite $\mathcal{H}^{n{-}2}$ measure).
In this way we separate the $C^1$ hypersurfaces one from each other.
Now $J_{u_k}$ is included in a finite union of \emph{pairwise disjoint} compact $C^1$ hypersurfaces contained in $\Omega$. It is then enough to apply \cite[Lemma~4.3]{DPFusPra17} to get a slight modification of $u_k$ such that $J_{u_k}$ indeed coincides with the finite union of $C^1$ hypersurfaces above.
Therefore the proof is concluded.
\end{proof}
We now start the proof of Theorem~\ref{teo:TeorB}. The following Lemma is employed in Proposition~\ref{prop:roughTeorB}, which is the counterpart of Theorem~\ref{teo:rough} in the proof of Theorem~\ref{teo:TeorB}.
\begin{lemma}\label{le:controllo eu}
Let $Q=(-2r,2r)^n$, $Q'=(-r,r)^n$, $v\in SBD^p_\infty(Q_r)$, and 
 $\varphi_r(x):=r^{-n}\varphi_1(x/r)$, with $\varphi_1 \in C_c^\infty(B_1)$. Then (recall that $\E^j v$ is the jump part of the measure $\E v$, see \eqref{1205181701})
\begin{equation}\label{1205180959}
\int \limits_{Q'} |e(v\ast \varphi_r) - e(v)\ast \varphi_r|^p \dx \leq   \|\varphi_1\|_{L^p(B_1)}^p\, r^{-n(p-1)}  \big(|\E^j v|(Q) \big)^p  \,.
\end{equation} 
\end{lemma}
\begin{proof}
From the standard approximation argument by Anzellotti and Giaquinta (cf.\ e.g.\ \cite[Theorem~5.2]{AnzGia82}) there exist $v_h\in C^\infty(Q;\Rn)\cap BD(Q)$ such that $v_k\to v$ in $L^1(Q;\Rn)$, there is the convergence in mass $\|e(v_k)\|_{L^1(Q)} \to |\E v|(Q)$, and
\begin{equation}\label{1205181643}
\|e(v_k-v)\|_{L^1(Q)} \to |\E^j v|(Q)\,.
\end{equation}
For any $k\in \N$ we have that
\begin{equation}\label{1205181655}
\|e(v_k\ast \varphi_r)-e(v)\ast \varphi_r\|_{L^p(Q')}= \|e(v_k-v) \ast \varphi_r\|_{L^p(Q')} \leq \|\varphi_r\|_{L^p(B_r)} \|e(v_k-v)\|_{L^1(Q)}\,.
\end{equation}
Moreover $v_k\ast \varphi_r \to v\ast \varphi_r$ uniformly in $Q$, since $v_k\to v$ in $L^1(Q;\Rn)$, and then \eqref{1205181655} implies that $e(v_k\ast \varphi_r)$ is bounded in $L^p$ with respect to $k$, so that
\begin{equation*}
e(v_k\ast \varphi_r)\weak e(v\ast \varphi_r) \quad \text{in }L^p(Q';\Mnn)\,.
\end{equation*}
We employ the convergence above to pass to the limit in the left hand side of \eqref{1205181655}, while for the right hand side we use \eqref{1205181643}, so
\begin{equation*}
\|e(v\ast \varphi_r)-e(v)\ast \varphi_r\|_{L^p(Q')} \leq \|\varphi_r\|_{L^p(B_r)} |\E^j v|(Q) \,.
\end{equation*}
Now \eqref{1205180959} follows raising to the $p$ and observing that
\begin{equation*}
\int \limits_{B_r} |\varphi_r|^p \dx = r^{-n p} \int \limits_{B_r} |\varphi_1(x/r)|^p \dx =  r^{-n(p-1)} \int \limits_{B_1} |\varphi_1|^p \,\mathrm{d}y\,.
\end{equation*}
\end{proof}
\begin{proposition}\label{prop:roughTeorB}
Let $\Omega$, $\widetilde{\Omega}$ be bounded open subsets of $\Rn$, with $\ol \Omega\subset \widetilde{\Omega}$, and let $u\in SBD^p_\infty(\widetilde{\Omega})$, $p>1$.
Then there exist $u_k\in SBV^p(\Omega;\Rn)\cap L^\infty(\Omega; \Rn)$ such that $J_{u_k}$ is included in a finite union of $(n-1)$--dimensional closed cubes, $u_k\in C^\infty(\ol\Omega\setminus J_{u_k}; \Rn) \cap W^{m,\infty}(\Omega\setminus J_{u_k}; \Rn)$ for every $m\in \N$, and:
\begin{subequations}
\begin{align}
\lim_{k\to \infty} \|u_k-& u\|_{L^1(\Omega;\Rn)}=0\,,\label{1roughTeorB}\\
\limsup_{k\to \infty} \int\limits_\Omega  |e(u_k)|^p \dx & \leq \int\limits_\Omega |e(u)|^p \dx + C |\E^j u|(\widetilde{\Omega})\,, \label{2roughTeorB}\\
\hn(J_{u_k})& \leq k |\E^j u|(\widetilde{\Omega})\,,\label{3roughTeorB}\\
\limsup_{k\to \infty} \int \limits_{J_{u_k}} \big| [u_k] \big| \dh  &\leq C \int \limits_{J_u} \big| [u] \big| \dh + C \int \limits_{\widehat{\Omega}_1} |e(u)|^p\dx\,, \label{4roughTeorB}
\end{align}
\end{subequations}
with $\mathcal{L}^n(\widehat{\Omega}_1)\leq C|\E^j u|(\widetilde{\Omega})$, and $C>0$ independent of $k$. 
\end{proposition}
\begin{proof}
As in Theorem~\ref{teo:rough}, we take $k\in \N$ with $k> \frac{8 \sqrt{n}}{\mathrm{dist }(\partial \Omega, \partial \widetilde{\Omega})}$, $\varphi\in C^\infty_c(B_1)$ radial, $\varphi_k(x)=k^n \varphi(kx)$, and for any $z\in (2 \km) \Z^n \cap \Omega$ the cubes \[q_z^k:=z+(-\km,\km)^n\,, \quad\tq_z^k:= z+(-2\km,2\km)^n\,.\]
We take the ``good'' and ``bad'' nodes
\begin{equation}\label{1205181848}
\widehat{G}^k:=\{z\in (2 \km) \Z^n \cap \Omega \colon |E^j u|(\tqz)\leq k^{-n}\}\,, \qquad\widehat{B}^k:= z\in (2 \km) \Z^n \cap \Omega \sm \widehat{G}^k\,,
\end{equation}
and the corresponding sets
\begin{equation}\label{1205181927}
\widehat{\Omega}^k_g:=\bigcup_{z\in \widehat{G}^k} \qz\,,\quad \widehat{\Omega}^k_b:= \bigcup_{z\in \widehat{B}^k} \tqz\,,
\end{equation}
so $\widehat{\Omega}^k_b= \widetilde{\Omega}\sm \widehat{\Omega}^k_g + (-k^{-1}, k^{-1})^n$.
We have (recall that $\tqz$ are finitely overlapping)
\begin{equation}\label{1205181928}
\# \widehat{B}^k\leq C |E^j u|(\widetilde{\Omega}) \, k^n\,,
\end{equation}
so that
\begin{equation*}
\mathcal{L}^n(\widehat{\Omega}^k_b)\leq C |E^j u|(\widetilde{\Omega})\,.
\end{equation*}
By Lemma~\ref{le:controllo eu} and \eqref{1205181848}, for every $z\in \widehat{G}^k$
\begin{equation}\label{1205181924}
\int \limits_{\qz} |e(u\ast \varphi_k) - e(u) \ast \varphi_k|^p\dx \leq C |\E^j u|(\tqz)\,.
\end{equation}
Notice that here this plays the same role of \eqref{prop3iiCCF16applicata} for Theorem~\ref{teo:density}.
We then define the approximating functions as
\begin{equation}\label{1205182113}
u_k:= \begin{cases}
u\ast \varphi_k \quad &\text{in }\Omega\sm \widehat{\Omega}^k_b\,,\\
\tilde{a}_z \quad &\text{in }\qz \cap \widehat{\Omega}^k_b\,,
\end{cases}
\end{equation}
where $\tilde{a}_z\colon \Rn \to \Rn$ is affine with $e(\tilde{a}_z)=0$ such that (cf.\ \eqref{eq:PoinBad})
\begin{equation*}
\| u-\tilde{a}_z \|_{L^1(\tilde{q}_z)} \leq C k^{-1} | \E u|(\tqz)\,.
\end{equation*}
It is not difficult to see that 
\begin{equation}\label{1205182131}
\|u-u_k\|_{L^1(\Omega;\Rn)}\leq C k^{-1} |\E u|(\widetilde{\Omega})\,, 
\end{equation}
so \eqref{1roughTeorB} follows.

As done for \eqref{1204181419} and \eqref{1204181414}, we have that 
\[
J_{u_k}\subset \bigcup_{z\in \widehat{B}_k}(J_{u_k}\cap \ol \tqz)\quad\text{ and, for $z\in \widehat{B}^k$, }\quad \hn(J_{u_k} \cap \ol\tqz)\leq C \, k^{n-1}\,,
\]
therefore \eqref{3roughTeorB} follows from \eqref{1205181928}.
Similarly to \eqref{1304180952} it follows that
\begin{equation*}
\int \limits_{\interior{(\widehat{\Omega}_b^k)}} \big| [u_k] \big| \dh \leq C\,|\E u|(\widehat{\Omega}_b^k)\,,
\end{equation*}
while for every $x\in \partial\widehat{\Omega}^k_b \cap \qz$
\begin{equation*}
\big|[u_k]\big|(x)=|(u-\tilde{a}_z)\ast \varphi_k|(x)\leq C k^{n} \|u-\tilde{a}_z\|_{L^1(\tqz)}\leq C k^{n-1}|\E u|(\tqz)\,.
\end{equation*}
Integrating the above inequality we deduce
\begin{equation*}
\int \limits_{\partial\widehat{\Omega}^k_b \cap \qz} \big| [u_k] \big| \dh \leq C |\E u|(\tqz)
\end{equation*}
and, since the cubes $\tqz$ are finitely overlapping and $J_{u_k}\subset \widehat{\Omega}^k_b$,
\begin{equation*}
\int \limits_{J_{u_k}} \big| [u_k] \big| \dh  \leq C \int \limits_{J_u \cap \widetilde{\Omega}} \big| [u] \big| \dh + C \int \limits_{\widehat{\Omega}^k_{b,1}} |e(u)|\dx\,,
\end{equation*} 
where $\widehat{\Omega}^k_{b,1}:= \widehat{\Omega}^k_b + (-k^{-1}, k^{-1})^n$. This gives \eqref{4roughTeorB} with $\widehat{\Omega}_1=\widehat{\Omega}^k_b$, since $e(u)\in L^p(\Omega;\Mnn)$.

We prove \eqref{2roughTeorB} by summing up \eqref{1205181924} over $z\in \widehat{G}^k$ (we use again that $\tqz$ are finitely overlapping) and recalling that $e(u_k)=0$ in $\widehat{\Omega}^k_b$, see \eqref{1205182113}. This concludes the proof.
\end{proof}
\begin{proof}[Proof of Theorem~\ref{teo:TeorB}]
As in Theorem~\ref{teo:TeorA} we follow the proof of Theorem~\ref{teo:density}, replacing the definition of $(u_k)_\m$ in \eqref{2004180730} by
\begin{equation}\label{1205182203}
(u_k)_{\m}:= k\text{-th approximating function for }u_{\m} \text{ on }Q^-_{\m}, \text{ by Proposition~\ref{prop:roughTeorB}}\,,
\end{equation}
and the definition of $(u_k)_{B_0}$ with Proposition~\ref{prop:roughTeorB} in place of Theorem~\ref{teo:rough}.

By \eqref{1205182131} we have the following analogue of \eqref{2204180755}
\begin{equation}\label{1205182207}
\|(\widehat{u}_{Q^-})-(u_k)_{Q^-}\|_{L^1(Q^-)} \leq  \,C k^{-1} |\E u|\big((Q^-+(-t,t)^n)\sm \Gamma\big)\,, 
\end{equation}
that at the end gives (by definition of $u_k$ and by Proposition~\ref{prop:roughTeorB} in $B_0$) 
\begin{equation}\label{1205182344}
u_k\to u \quad\text{in }L^1(\Omega;\Rn)\,.
\end{equation}

On the other hand, \eqref{2roughTeorB} implies the following counterpart of \eqref{1105182102}
\begin{equation*}
\|e((u_k)_{Q^-})\|^p_{L^p(Q^-)} \leq \|e(\widehat{u}_{Q^-})\|^p_{L^p(Q^-+(-t,t)^n)} + C |\E^j u|(Q^-\sm \Gamma) \,,
\end{equation*}
and similarly in $B_0$. Therefore we deduce
\begin{equation}\label{1205182345}
\limsup_{k\to \infty} \|e(u_k)\|^p_{L^p(\Omega;\Mnn)} \leq \|e(u)\|^p_{L^p(\Omega;\Mnn)} + C |\E^j u|(\Omega\sm \widehat{\Gamma}) \,,
\end{equation}
and the last term goes to 0 as $\varepsilon\to 0$ by \eqref{1205182222}.

Let us now consider $[u_k]$. In comparison to \eqref{2104182323}, \eqref{4roughTeorB} gives also an additional term 
\begin{equation*}
C \int \limits_{\widehat{\Omega}^-_{1,\m}} |e(u)|^p\dx\,,
\end{equation*}
with $\mathcal{L}^n(\widehat{\Omega}^-_{1,\m})\leq C |\E^j u|(Q^-_\m)$. Summing up on $\m$ this entails in \eqref{2204181346} an additional term
\begin{equation*}
C \int \limits_{\widehat{\Omega}_1} |e(u)|^p\dx\,,
\end{equation*}
with $\mathcal{L}^n(\widehat{\Omega}_1)\leq C |\E^j u|(\Omega\sm \widehat{\Gamma})$, that goes to 0 in $\varepsilon$ by \eqref{1205182222}.

The estimate of $[u_k]$ on $J'_\m$ is done as in \eqref{2604182001}, distinguishing four cases according to the fact that each cube intersecting $J'_\m$ is good or bad with respect to $Q^-_\m$ or $Q^-_{\m+e_i}$. The difference is in the definition of $(u_k)_{\m}$: now there are no exceptional sets in the good cubes, but $[u]$ enters also if a cube is good regarded both in $Q^-_\m$ and in $Q^-_{\m+e_i}$ (we employ \eqref{1205181924} in place of \eqref{prop3iiCCF16applicata}). The final estimate is anyway the same of \eqref{2604182001}, and this holds also for \eqref{2204180925}. Then we obtain as in \eqref{2204181346} that (take always $\varepsilon\to 0$ more slowly than $k^{-1}$)
\begin{equation*}
\lim_{k\to \infty} \int \limits_{J_{u_k}\sm \widehat{\Gamma}} \big| [u_k] \big| \dh =0\,.
\end{equation*}

In the same way also the estimate \eqref{2204180800} is still true, and combined with \eqref{1205182207} this implies
\begin{equation*}
\begin{split}
\int \limits_{\widehat{\Gamma}} \big|[u]-[u_k]\big| \dh  \leq C\, \|e(u)\|_{L^1(\{\mathrm{d}(\cdot, \widehat{\Gamma})<t\})} + C\, |\E^j u| (\Omega\sm \widehat{\Gamma}) + \frac{C}{t} k^{-1} |\E u| (\Omega\sm \widehat{\Gamma}) \,.
\end{split}
\end{equation*}

Then, in particular, $|E^j (u-u_k)|(\Omega) \to 0$, and \eqref{1205182344}, \eqref{1205182345} give $u_k$ bounded in $BD(\Omega)$ and thus \eqref{eqs:mainTeorB} by \cite[Theorem~1.1]{BelCosDM98}. The proof is then concluded.
\end{proof}

\begin{remark}\label{rem:1805180933}
As in \cite[Theorem~B]{DPFusPra17}, that deals with $SBV^p_\infty$ functions, we are not able to ensure that $\hn(J_{u_k}\sm J_u) \to 0$ in Theorem~\ref{teo:TeorB}. This comes from \eqref{3roughTeorB}, which in turn is a consequence of \eqref{1205180959} in Lemma~\ref{le:controllo eu}. Improving this estimate could then give a control on the measure of the jump created in the approximation procedure.
\end{remark}

\begin{remark}\label{rem:1705181900}
In Theorems~\ref{teo:density} and \ref{teo:TeorB} the jump of the approximating functions is contained in a finite union of $C^1$ hypersurfaces, which are not necessarily pairwise disjoint. Indeed, an issue comes from the intersections of $\widehat{\Gamma}$ with the bad (and the boundary good) cubes coming from the construction in Theorem~\ref{teo:rough} and Proposition~\ref{prop:roughTeorB} in any $Q^-_\m$:
if it is possible to choose the cubes of sidelength $k^{-1}$ in such a way that the grid intersects $\widehat{\Gamma}$ (and $\widehat{\Gamma}_{\partial\Omega}$) in a finite number of pairwise disjoint components of finite $\mathcal{H}^{n{-}2}$-measure (this should be guaranteed by a delicate use of the area formula for Lipschitz graph, since $\widehat{\Gamma}$ is a finite union of pairwise disjoint $C^1$ curves), then one could use the capacitary argument in \cite[Corollary~3.11]{Cor97} if $p\in (1,2]$ to replace the jump on this $(n{-}2)$-dimensional set by a smooth transition, so separating the hypersurfaces. For $p>2$ the situation is more delicate since one can apply \cite[Lemma~5.2]{DPFusPra17} only if $J_u\subset \subset \Omega$ and $u\in C^1(\Omega\sm \ol{J_u})$. On the other hand, one could argue as in Theorem~C of \cite{DPFusPra17}, in Part~B--Steps II, III (see Remark~\ref{rem:1705181845}) to separate $J_u$ from $\partial\Omega$, but losing $u\in C^1$ near $\partial \Omega$. Here we choose to avoid this possible refinement due to these technicalities and since in the applications considered (also in \cite{DPFusPra17}) one needs just $J_u$ closed or one passes through the approximation in \cite{CorToa99}, that permits to separate the components.
\end{remark}

\section{Some applications}\label{sec:Appl}


The theorems of this paper on $SBD$ functions may be employed in combination with other density result in $SBV$, such as those 
in \cite{BraChP96}, \cite{CorToa99}, or \cite{DPFusPra17}. In particular, 
Cortesani and Toader approximate functions in $SBV^p(\Omega;\Rn) \cap L^\infty(\Omega;\Rn)$ by so-called ``piecewise smooth''  $SBV$-functions, denoted $\mathcal{W}(\Omega;\Rn)$, namely
\begin{equation*}
u \in \mathcal{W}(\Omega;\Rn) \text{ if }\begin{cases}
u\in SBV(\Omega;\Rn)\cap W^{m,\infty}(\Omega\sm J_u;\Rn) \,\text{for every }m\in \N\,,\\
\hn(\ol{J_u} \sm J_u ) = 0\,,\\
\ol{J_u} \text{ is the intersection of $\Omega$ with a finite union of ${(n{-}1)}$-dimensional simplexes}\,.
\end{cases}
\end{equation*}
We report below the result by Cortesani and Toader, in a slightly less general version.
\begin{theorem}[\cite{CorToa99}, Theorem~3.1] \label{teo:CorToa}
Let $\Omega$ be an open bounded Lipschitz set.
For every $u\in SBV^p(\Omega;\Rn) \cap L^\infty(\Omega;\Rn)$ there exist $u_k\in \mathcal{W}(\Omega;\Rn)$ such that
\begin{align*}
\lim_{k\to \infty} \Big( \|u_k-u\|_{L^1(\Omega;\Rn)} &+ \|\nabla u_k -\nabla u\|_{L^p(\Omega;\mathbb{M}^{n\times n})} + \hn(J_{u_k}\triangle J_u)  \Big)=0\,,\\
\lim_{k \to \infty} \int \limits_{J_{u_k}\cap A} \phi(x, u_k^+, & u_k^-, \nu_{u_k})   \dh = \int \limits_{J_u \cap A} \phi(x, u^+, u^-, \nu_{u}) \dh\,,
\end{align*}
for every $A\subset \Omega$, $\hn(\partial A \cap J_u)=0$, and every $\phi$ strictly positive, continuous, and $BV$-elliptic (see e.g.\ \cite{Amb90GSBV} or \cite[equation (2.4)]{CorToa99} for the notion of $BV$-ellipticity).
\end{theorem}
\begin{remark}\label{rem:1705181845}
During the proof of Theorem~C of \cite{DPFusPra17}, in Part~B--Steps II, III, it is shown that for every $\varepsilon>0$ and $u \in SBV^p(\Omega;\Rn) \cap W^{1,\infty}(\Omega\sm J_u;\Rn)$ with $J_u$ closed, there is a $v$ with the same regularity, such that $J_v \subset \subset \Omega$ and $\big(\|u-v\|_{BV} + \|\nabla (u-v)\|_{L^p} + \hn(J_u \triangle J_v)\big)<\varepsilon$. Moreover, by the  procedure of \cite[Theorem~3.1]{CorToa99}, the function $v$ may be approximated in the sense of Theorem~\ref{teo:CorToa} by $v_k\in\mathcal{W}(\Omega;\Rn)$ such that also $J_{v_k} \subset \subset \Omega$. Then by a diagonal argument we may assume that $J_{u_k}\subset \subset \Omega$ in Theorem~\ref{teo:CorToa}.
\end{remark}
Theorems~\ref{teo:density} and \ref{teo:CorToa} are, in particular, very useful tools to prove $\Gamma$-convergence approximations for energies including a bulk part depending on $e(u)$ and a surface part depending on the measure of the jump set and on the amplitude of the jump. These energies are then formulated in the space $SBD^p$ and arise in particular in Fracture Mechanics. Indeed, the jump set may represent the set where a material is cracked, so that the surface part is usually interpreted as a dissipative part. In the present context we consider the case where the dissipation actually depends on the amplitude of the jump.
If the dissipation depends only on the measure of the jump set the fracture is said ``brittle'', in the other cases it is often called ``cohesive''.

The use of Theorems~\ref{teo:density} and \ref{teo:CorToa} permits to prove the $\Gamma$-limsup inequality just for $\mathcal{W}(\Omega;\Rn)$ functions: one may approximate any $u \in SBD^p$ by $\widehat{u}_k\in \mathcal{W}(\Omega;\Rn)$, and, if one knows how to construct a \emph{recovery sequence} for functions in $\mathcal{W}(\Omega;\Rn)$, a diagonal argument is sufficient to conclude.

As an application of this strategy, we extend the following two results, for which the corresponding $\Gamma$-limsup inequality is proven in $\mathcal{W}(\Omega;\Rn)$ (and then extended to $SBD^p(\Omega) \cap L^\infty(\Omega;\Rn)$ by \cite{Iur14}). We notice that when the bulk energy depends on $e(u)$ it is not natural to assume that the minimisers are bounded, even if the boundary datum is bounded. Indeed, the functional is not only non decreasing by truncation, but it is not even true that a truncation of a $BD$ function is still in $BD$.

The first result is shown by Focardi and Iurlano in \cite[Theorem~3.2]{FocIur14}. Its generalisation is the following. (We formulate the result in a slightly less general setting to simplify the notation.)
\begin{theorem}\label{teo:genFI14}
Let $\Omega$ be an open bounded Lipschitz set, let $p>1$, $p':=p/(p-1)$, and $\psi \in C([0,1])$ decreasing with $\psi(1)=0$. Then the functionals $F_\varepsilon\colon L^1(\Omega;\Rn){\times}L^1(\Omega)$ defined as
\begin{equation*}
F_\varepsilon(u,v):=\begin{dcases}
\int \limits_\Omega \Big( v\, |e(u)|^2 +\frac{\psi(v)}{\varepsilon}+\varepsilon^{p-1} |\nabla v|^p \Big) \dx \quad &\text{if }(u,v)\in H^1(\Omega;\Rn){\times}W^{1,p}(\Omega;[\varepsilon,1])\,,\\
+\infty &\text{otherwise,} 
\end{dcases}
\end{equation*}
$\Gamma$-converge, as $\varepsilon\to 0$, in $L^1(\Omega;\Rn){\times}L^1(\Omega)$ to
\begin{equation*}
F(u,v):= \begin{dcases}
\int \limits_\Omega |e(u)|^2 \dx + a \hn(J_u) + b \int \limits_{J_u} \big| [u] \odot \nu_u \big| \dh \quad&\text{if }u\in SBD^2(\Omega),\, v=1\,,\\
+\infty &\text{otherwise,}
\end{dcases}
\end{equation*}
where 
$a:= 2 p^{1/p} p'^{1/p'} \int_0^1 \psi^{1/p}(s) \,\mathrm{d}s$ and $b:= 2 \psi^{1/2}(0)\,.$
\end{theorem} 
\begin{remark}
In \cite[Remark~4.5]{FocIur14} the authors explain why it was possible to prove the $\Gamma$-limsup inequality only with an \emph{a priori} $L^\infty$ bound on $u$. Here we improve also the desired result in \cite[Remark~4.5]{FocIur14}, since we not only show it for $u\in SBD^2(\Omega)\cap L^2(\Omega;\Rn)$, but directly in $SBD^2(\Omega)$, without any additional integrability assumption. Notice that Theorem~\ref{teo:rough} would give a density result in $SBD^p(\Omega) \cap L^p(\Omega;\Rn)$ with the approximation technique in \cite{Cha04, Iur14} based on gluing rough approximations by means of a partition of unity. The work done in Section~\ref{sec:proofMain} is devoted to remove even the \emph{a priori} $L^p$ bound.
\end{remark}

We consider now a result proven very recently by Caroccia and Van Goethem, that enrichs \cite[Theorem~3.2]{FocIur14} with the presence of a low order potential $P$, controlled from above and below by two linear functionals in $e(u)$. This is related to the simulation of models for fluid-driven fracture (e.g.\ \emph{fracking} and hydraulic fracture in porous media), and goes in the direction of the treatment of non-interpenetration or Tresca-type conditions for plastic slips. The result is \cite[Theorem~2.3]{CarVG18}, and the $\Gamma$-limsup inequality is still proven for $u\in \mathcal{W}(\Omega;\Rn)$.
We state below directly the generalised result simplifying some notation, as done for Theorem~\ref{teo:genFI14}.
\begin{theorem}\label{teo:genCarVG18}
Let $\Omega$ be an open bounded Lipschitz set, $\psi \in C([0,1])$ decreasing with $\psi(1)=0$, $P\colon \Omega{\times}\Mnn\to \R$ continuous in the first argument, convex in the second, with $-\sigma|M|\leq P(x,M) \leq l|M|$ for any $l>0$ and a suitable $0<\sigma< 2 \psi^{1/2}(0)$. Then the functionals $G_\varepsilon\colon L^1(\Omega;\Rn){\times}L^1(\Omega)$ defined as
\begin{equation*}
G_\varepsilon(u,v):=\begin{dcases}
\int \limits_\Omega \Big( v\, |e(u)|^2 +\frac{\psi(v)}{\varepsilon}+P(x, e(u)) \Big) \dx \quad &\text{if }(u,v)\in H^1(\Omega;\Rn){\times}V_\varepsilon\,,\\
+\infty &\text{otherwise,} 
\end{dcases}
\end{equation*}
where
\begin{equation*}
V_\varepsilon:=\{v\in W^{1,\infty}(\Omega; [\varepsilon,1]) \colon |\nabla v|\leq 1/\varepsilon\}\,,
\end{equation*}
$\Gamma$-converge, as $\varepsilon\to 0$, in $L^1(\Omega;\Rn){\times}L^1(\Omega)$ to $G(u,v)$ given by
\begin{equation*}
\begin{dcases}
\int \limits_\Omega \hspace{-0.5em} \Big(|e(u)|^2 + P(x, e(u))\Big) \dx +  \int \limits_{J_u} \hspace{-0.5em}\Big( a' + b' \big| [u] \odot \nu_u \big| + P_\infty(\cdot, [u] \odot \nu_u) \Big)\dh &\text{in }SBD^2(\Omega){\times}\{v=1\},\\
+\infty &\text{otherwise,}
\end{dcases}
\end{equation*}
where $a':= 2 \int_0^1 \psi(s) \, \mathrm{d}s$, $b':=2 \psi^{1/2}(0)$, and
$P_\infty(x,M):=\lim_{t\to +\infty} \frac{P(x,t M)-P(x,0)}{t}$.
\end{theorem}
We conclude noticing that with our result it is possible to deal with bulk energies having growth $p>1$ in $e(u)$, and not necessarily quadratic.
As observed in \cite{CFI17Density}, the constructions by \cite{Cha04} and \cite{Iur14} do not provide approximations in $(G)SBD^p$ but only in $(G)SBD^2$.
From a mechanical point of view the $p$-growth of the bulk energy is connected with elasto-plastic materials (see for instance \cite[Sections~10 and 11]{Hut} and reference therein) and interesting also in a purely elastic framework (see \cite[Section~2]{CFI17Density}).

In a future paper 
functionals with non quadratic bulk energy and dissipated energy depending only on the deviatoric part of the matrix-valued function $[u]\odot \nu_u$ will be investigated.

%
%
%
%
%


\bigskip
\noindent {\bf Acknowledgements.}
I have been supported by a public grant as part of the \emph{Investissement d'avenir} project, reference ANR-11-LABX-0056-LMH, LabEx LMH, and acknowledge the financial support from the Laboratory Ypatia and the CMAP.
I am grateful to Antonin Chambolle for his generous and fruitful advices. 

%

\end{document}